\documentclass[square]{imsart}

\RequirePackage[OT1]{fontenc}
\RequirePackage{natbib}
\RequirePackage[colorlinks,citecolor=blue,urlcolor=blue]{hyperref}
\RequirePackage{hypernat}

\makeatletter \@addtoreset{equation}{section} \makeatother
\usepackage{graphicx}
\usepackage{amsmath}
\usepackage{amssymb}
\usepackage{amsfonts}
\usepackage{amsthm}
\usepackage{thmtools}
\usepackage{caption}
\usepackage{subcaption}

\usepackage{cleveref}
\usepackage[noend]{algpseudocode}
\usepackage{algorithm}
\usepackage[leftcaption]{sidecap}

\newcommand{\RR}{{\mathbb  R}}

\newtheorem{lem}{Lemma}[section]

\newtheorem{prop}{Proposition}[section]
\newtheorem{exmpl}{Example}[section]

\newtheorem{rem}{Remark}[section]

\newtheorem{cor}{Corollary}[section]

\advance \topmargin by -\headheight
\advance \topmargin by -\headsep
\advance \topmargin .2in
\begin{document}

\begin{frontmatter}
\title{The Bennett-Orlicz norm}  
\runtitle{Bennett-Orlicz norm} 

\begin{aug}
\author{\fnms{Jon A.} \snm{Wellner}\thanksref{t1}\ead[label=e1]{jaw@stat.washington.edu}}
\ead[label=u1,url]{http://www.stat.washington.edu/jaw/}

\thankstext{t1}{Supported in part by NSF Grants DMS-1104832 and DMS 1566514, and NI-AID grant 2R01 AI291968-04} 
\runauthor{Wellner}

\address{Department of Statistics \\University of Washington\\Seattle, WA 98195-4322\\
\printead{e1}}

\printead{u1}
\end{aug}

\begin{abstract}
Lederer and van de Geer (2013) introduced a new Orlicz norm, the Bernstein-Orlicz norm, which 
is connected to Bernstein type inequalities.
Here we introduce another Orlicz norm, the Bennett-Orlicz norm, which is connected to 
Bennett type inequalities.   
The new Bennett-Orlicz norm yields  inequalities for expectations of maxima which are potentially 
somewhat tighter than those resulting from the Bernstein-Orlicz norm when they are both applicable.
We discuss cross connections between these norms, exponential inequalities
of the Bernstein, Bennett, and Prokhorov types, and make comparisons with results of 
Talagrand (1989, 1994), and 
Boucheron, Lugosi, and Massart (2013).  

\end{abstract}



\end{frontmatter}

\tableofcontents

\bigskip

\bigskip

\section{Orlicz norms and maximal inequalities}
\label{sec:intro}
Let $\Psi $ be an increasing convex function from $[0,\infty)$ onto $[0,\infty)$.
Such a function is called a {\sl Young-Orlicz modulus} by \cite{MR1720712}, and a {\sl Young modulus} by \cite{MR1666908}.  
Let $X$ be a random variable.
The {\sl Orlicz norm} $\| X \|_{\Psi}$ is defined by 
\begin{eqnarray*}
\| X \|_{\Psi} = \inf \left \{ c>0 : \ E \Psi \left ( \frac{|X|}{c} \right ) \le 1 \right \} ,
\end{eqnarray*}
where the infimum over the empty set is $\infty$.  By Jensen's inequality it is easily shown that 
this does define a norm on the set of random variables for which $\| X \|_{\Psi}$ is finite.  
The most important functions $\Psi$ for a variety of applications are those of the form
$\Psi (x) = \exp ( x^p ) -1 \equiv \Psi_p (x)$ for $p\ge 1$, and in particular $\Psi_1$ and $\Psi_2$
corresponding to random variables which are ``sub-exponential'' or ``sub-Gaussian'' respectively.   
See 
\cite{MR0126722},  
\cite{MR1720712}, 
\cite{MR1348376}, 
\cite{MR1666908},   
and \cite{MR1385671}  
for further background on Orlicz norms, and see 
\cite{MR1113700},   
\cite{MR0126722},    
and \cite{MR0367121}  
for more information about Birnbaum-Orlicz spaces.

The following useful lemmas are from 
\cite{MR1385671}, pages 95-97, and \cite{MR1348376} 
(see also \cite{MR1666908}, pages 188-190), respectively.
\smallskip

\begin{lem}
\label{lem:GeneralOrliczMaxIneqFiniteList}
Let $\Psi$ be a convex, nondecreasing, nonzero function with 
$\Psi (0) = 0$ and $\limsup_{x,y\rightarrow \infty} \Psi (x) \Psi (y) / \Psi (c x y) < \infty$ for some 
constant $c$.  Then, for any random variables $X_1, \ldots , X_m$,
\begin{eqnarray}
\big \| \max_{1 \le j \le m} X_i \big \|_{\Psi} \le K \Psi^{-1} (m) \max_{1\le j \le m} \| X_i \|_{\Psi} 
\label{GenericOrliczMaximalInequalFinite}
\end{eqnarray}
where $K$ is a constant depending only on $\Psi$. \\
\end{lem}

\begin{lem}
\label{lem:ExpTypeOrliczMaxIneqSequence}
Let $\Psi$ be a Young Modulus satisfying 
\begin{eqnarray*}
\limsup_{x,y\rightarrow \infty} \frac{\Psi^{-1} (xy)}{\Psi^{-1} (x) \Psi^{-1} (y)} < \infty  \qquad \mbox{and} \qquad
\limsup_{x\rightarrow \infty} \frac{\Psi^{-1} (x^2)}{\Psi^{-1} (x)} < \infty .
\end{eqnarray*}
Then for some constant $M$ depending only on $\Psi$ and every sequence of random variables $\{ X_k : \ k \ge 1 \}$,
\begin{eqnarray}
\bigg \| \sup_{k\ge 1} \frac{| X_k |}{\Psi^{-1} (k)} \bigg \|_{\Psi} \le M \sup_{k \ge 1} \| X_k \|_{\Psi} .
\label{ExpTypeOrliczMaximalIneqSeq}
\end{eqnarray}
\end{lem}

The inequality (\ref{GenericOrliczMaximalInequalFinite}) 
shows that if Orlicz norms for individual random  variables $\{ X_i \}_{i=1}^m $ 
are under control, then the $\Psi-$Orlicz norm of the maximum of the $X_i$'s is controlled by
a constant times $\Psi^{-1} (m)$ times the maximum of the individual Orlicz norms.  
The inequality (\ref{ExpTypeOrliczMaximalIneqSeq}) shows a stronger related Orlicz norm control of the 
supremum of an entire sequence $X_k$ divided by $\Psi^{-1} (k)$ if the supremum of the individual Orlicz 
norms is finite.   Lemma~\ref{lem:ExpTypeOrliczMaxIneqSequence} implies Lemma~\ref{lem:GeneralOrliczMaxIneqFiniteList} 
for Young functions of exponential type (such as $\Psi_p (x) = \exp (x^p) -1$ with $p \ge 1$), but it does not hold 
for power type Young functions such as $\Psi (x) = x^p$, $p \ge 1$.   These latter Young functions continue to be covered 
by Lemma~\ref{lem:GeneralOrliczMaxIneqFiniteList}.  
\cite{MR1348376}
carefully define Young moduli $\Psi_p (x) = \exp (x^p) -1$  for all $p>0$ and use 
Lemma~\ref{lem:ExpTypeOrliczMaxIneqSequence} to establish laws of the iterated logarithm for U-statistics.

A general theme is that  if $\Psi_a \le \Psi_b$
and we have control of the individual $\Psi_b$ Orlicz norms, then 
Lemma~\ref{lem:GeneralOrliczMaxIneqFiniteList} 
or Lemma~\ref{lem:ExpTypeOrliczMaxIneqSequence} applied with $\Psi = \Psi_b$ 
will yield a better bound than with $\Psi = \Psi_a$ in the sense that $\Psi_b^{-1} (m) \le \Psi_a^{-1} (m)$.

Here we are interested in functions $\Psi $ of the form 
\begin{eqnarray}
\Psi (x) = \exp ( h(x) ) -1
\label{OrliczFunctionNotAPower}
\end{eqnarray}
where $h$ is a nondecreasing convex function with $h(0) = 0$ {\sl not of the form} $x^p$.  
In fact, the particular functions $h$ of interest here are (scaled versions of):
\begin{eqnarray*}
h_0(x) & = & \frac{x^2}{2(1+x)}, \\
h_1 (x) & = & 1+x - \sqrt{1+2x} , \\
h_2 (x) & = & h(1+x) = (1+x) \log (1+x) - x, \\
h_4 (x) &  = & (x/2) \mbox{arcsinh} (x/2) , \\
h_5 (x) & = & (x) \mbox{arcsinh}(x/2) - 2 \Big ( \mbox{cosh} \big ( \mbox{arcsinh} (x/2)\big ) -1 \Big ) 
\end{eqnarray*}
for the particular $h(x) \equiv x (\log x -1 )+1$.    The functions $h_0$ and $h_1$ are related to 
Bernstein exponential bounds and refinements thereof due to \cite{MR1653272},   
while the function $h_2$ is related to Bennett's inequality  (\cite{Bennett:62}), and $h_4$ is related to Prokhorov's inequality
(\cite{MR0121857}).

\cite{MR3101846} 
studied the family of Orlicz norms defined in terms of scaled versions of $h_1$, and called 
called  them {\sl Bernstein-Orlicz norms}.    Our primary goal here is to compare and contrast the 
Orlicz norms defined in terms of $h_0$, $h_1$, $h_2$, and $h_4$. 
We begin in the next section 
by reviewing the Bernstein-Orlicz norm(s) as defined by \cite{MR3101846}.   
Section~3  gives corresponding results for what we call the {\sl Bennett-Orlicz norm(s)} corresponding 
to the function $h_2$.  
In Section~4 we give further comparisons and two applications.  

\section{The Bernstein-Orlicz norm}
\label{BernsteinOrlicz}
For a given number $L>0$, \cite{MR3101846}  
have defined the Bernstein-Orlicz norm  $\| X \|_{\Psi_L}$ with 
\begin{eqnarray}
\Psi_L (x) \equiv \Psi_1 (x; L) \equiv \exp \left \{ \left ( \frac{\sqrt{1+2 Lx}-1}{L} \right )^2 \right \} -1 
= \exp \left \{ \frac{2}{L^2} h_1 (Lx) \right \} -1 .
\label{PsiFcn-BernsteinOrlicz}
\end{eqnarray}
It is easily seen that 
\begin{eqnarray*}
\Psi_1(x;L) \sim \left \{ \begin{array}{l l} \exp (x^2) -1  & \ \ \mbox{for} \ Lx \ \mbox{small}, \\
                                                            \exp (2 x /L) -1 & \ \ \mbox{for} \ Lx \ \mbox{large} .
                                  \end{array} \right . 
\end{eqnarray*}

The following three lemmas of \cite{MR3101846} 
should be compared with the development on
page 96 of \cite{MR1385671}.  
\medskip

\begin{lem}
\label{lem:LvdG-lem1} 
Let $\tau \equiv \| Z \|_{\Psi_1 (\cdot; L)}$.  Then
\begin{eqnarray*}
P( | Z | > \tau [ \sqrt{t} + 2^{-1} L t ]) \le 2 e^{-t} \ \ \mbox{for all} \ t>0;
\end{eqnarray*}
or, equivalently, with $h_1^{-1} (y) \equiv y + \sqrt{2y}$, 
\begin{eqnarray*}
P \left  ( | Z | > (\tau/L) h_1^{-1} \left ( \frac{L^2 t}{2} \right ) \right ) \le 2 e^{-t} \ \ \mbox{for all} \ t>0;
\end{eqnarray*} 
or
\begin{eqnarray}
P( | Z | > z ) \le 2 \exp \left ( - \frac{2}{L^2} h_1 \left ( \frac{Lz}{\tau} \right ) \right )  \ \ \mbox{for all} \ z > 0 .
\label{H1_BernsteinInequal}
\end{eqnarray}
\end{lem}
\medskip
  
\begin{lem} 
\label{lem:LvdG-lem2}
Suppose that for some $\tau $ and $L>0$ we have 
\begin{eqnarray*}
P( | Z | \ge \tau [ \sqrt{t} + 2^{-1} L t]) \le 2 e^{-t}   \ \ \mbox{for all} \ t>0 .
\end{eqnarray*}
Equivalently,  the inequality (\ref{H1_BernsteinInequal}) holds.  
Then $ \| Z \|_{\Psi_1(\cdot ;\sqrt{3}L)} \le \sqrt{3} \tau $. 
\end{lem}
\smallskip

\begin{exmpl}
\label{ex:PoissonBernOrlicz}
Suppose that $X \sim \mbox{Poisson} (\nu)$.  Then it is well-known
(see e.g. \cite{MR3185193}, page 23), 
that 
\begin{eqnarray*}
P( | X - \nu | \ge z ) \le 2 \exp \left ( - \nu h_2 (z / \nu) \right ) \le 2 \exp \left (- 9 \nu h_1 (z/(3 \nu)) \right ) 
\end{eqnarray*}
where $h_2 (x) = h(1+x) = (x+1) \log (x+1) - x$.  
Thus the inequality involving $h_1$ holds with $9 \nu = 2/L^2$ and $1/(3\nu ) =  L / \tau$.
Thus $L = \sqrt{2/(9 \nu)} = 3^{-1} \sqrt{2/\nu}$, $\tau = L 3 \nu = \sqrt{2\nu}$.  We conclude from 
Lemma~\ref{lem:LvdG-lem2} that 
\begin{eqnarray*}
\| X - \nu \|_{\Psi_1 (\cdot; \sqrt{ 2/ 3\nu})} \le \sqrt{6 \nu} .
\end{eqnarray*}
\end{exmpl}

\cite{MR717231}  
and 
\cite{MR1089429} 
showed how to 
bound the Orlicz norm of the maximum of random variables with bounded Orlicz norms; 
see also 
\cite{MR1666908},  section 4.3, 
and 
\cite{MR1385671}, Lemma 2.2.2, page 96.    
The following bound for the expectation of the maximum 
was given by \cite{MR3101846}; also see \cite{MR3185193}, Theorem 2.5, pages 32-33.
\medskip

\begin{lem}
\label{lem:LvdG-lem3} 
Let $\tau $ and $L$ be positive constants, and let $Z_1, \ldots , Z_m$ 
be random variables satisfying $\max_{1 \le j \le m} \| Z_j \|_{\Psi_L} \le \tau$. 
Then
\begin{eqnarray}
E \{ \max_{1 \le j \le m} | Z_j | \} \le \tau \Psi_1^{-1} (m;L) = \tau \left \{ \sqrt{\log (1+m)} + \frac{L}{2} \log (1+m) \right \}.
\label{BernsteinOrliczMaxIneq}
\end{eqnarray}
\end{lem}
\medskip

\begin{cor}
\label{cor:JustLogTermBernsteinOrliczNormMax}
For $m\ge 2$
\begin{eqnarray*}
E \{ \max_{1 \le j \le m} | Z_j | \} \le \tau \Psi_1^{-1} (m;L) 
\le 2 \max \{ \tau , L \tau /2 \} \log (1 + m) .
\end{eqnarray*}
In particular when $Z_j \sim \mbox{Poisson}(\nu)$ for $1 \le j \le m$ 
\begin{eqnarray*}
E \{ \max_{1 \le j \le m} | Z_j | \} \le \tau \Psi_1^{-1} (m;L)  \le 2 \{ \sqrt{2 \nu} \vee 1/3 \} \log (1+m).
\end{eqnarray*}
\end{cor}

\par\noindent
{\bf Proof.}  This follows from Lemma~\ref{lem:LvdG-lem3} since $\sqrt{x} \le x$ for $x \ge 1$.   
The Poisson$(\nu)$ special case then follows from Example~\ref{ex:PoissonBernOrlicz}.
\hfill $\Box$
\medskip

It will be helpful to relate $\Psi_1(\cdot; L)$ to several functions appearing frequently 
in the theory of exponential bounds as follows:
for $x\ge 0$, we define
\begin{eqnarray}
&& h(x) = x(\log x -1 ) +1, \nonumber \\
&& h_1 (x) = 1+x - \sqrt{1+2x},  \label{hFcns}\\
&& h_0 (x) = \frac{x^2}{2(1+x)} .\nonumber
\end{eqnarray}
It is easily shown (see e.g. \cite{MR3185193} Exercise 2.8, page 47) that 
\begin{eqnarray}
9 h_0 (x/3) \le 9 h_1 (x/3) \le h(1+x) .
\label{hFcnInequalitiesOne}
\end{eqnarray}
A trivial restatement of the inequality on the left above and some algebra and easy inequalities 
yield
\begin{eqnarray}
h_0 (x) \le h_1 (x) \le 2 h_0 (x) \le h_0 (2x) .
\label{hFcnInequalitiesTwo}
\end{eqnarray}
The latter inequalities imply that the Orlicz norms based on $h_0$ and $h_1$ are equivalent up to constants.

One reason the functions $h_0$ and $h_1$ are so useful is that they both have explicit inverses:
from Boucheron, Lugosi, and Massart (2013), page 29, for $h_1$ and direct 
calculation for $h_0$,
\begin{eqnarray*}
&& h_1^{-1} (y) = y + \sqrt{2y} ,\ \ \mbox{for} \ \ y\ge 0,   \label{InverseOFhOne}\\
&& h_0^{-1} (y) = y+ \sqrt{y^2 + 2y}   \label{InverseOFhZero}.
\end{eqnarray*}

To relate the inequalities in Lemmas~\ref{lem:LvdG-lem1} and \ref{lem:LvdG-lem2} 
to more standard inequalities 
(with names) we note that 
\begin{eqnarray*}
\sqrt{t} + 2^{-1} L t = \frac{1}{L} h_1^{-1} \left ( \frac{L^2 t}{2} \right ) .
\end{eqnarray*}
This implies immediately that the inequality in 
Lemma~\ref{lem:LvdG-lem2} can be rewritten as
\begin{eqnarray*}
P( | Z | > z) 
& \le & 2 \exp \left ( - \frac{2}{L^2} h_1 \left ( \frac{Lz}{\tau} \right ) \right )  \\
& \le & 2 \exp \left ( - \frac{2}{L^2} h_0 \left ( \frac{Lz}{\tau} \right ) \right )  
 =  2 \exp \left ( - \frac{z^2}{\tau^2 + L \tau z } \right ) \ \ \ \ \mbox{for all} \ z>0.
\end{eqnarray*}

Here is a formal statement of a proposition relating exponential tail bounds in the traditional 
Bernstein form in terms of $h_0$ to tail bounds in terms of the (larger) function $h_1$. 
\medskip

\par\noindent
\begin{prop}
\label{prop:hOne-hZeroEquiv}
Suppose that a random variable $Z$ satisfies 
\begin{eqnarray}
P( | Z | > z) \le 2 \exp \left ( - \frac{z^2}{2 (A + Bz)} \right ) = 2 \exp \left ( - \frac{A}{B^2} h_0 \left ( \frac{Bz}{A} \right ) \right ) 
\ \ \ \mbox{for all} \ \ z > 0
\label{BernsteinTypeInequality}
\end{eqnarray}
for numbers $A,B>0$.  Then the hypothesis of 
Lemma~\ref{lem:LvdG-lem2} holds with $L$ and $\tau$ given by 
$L^2 = 2B^2/A$ and $\tau = 2^{3/2} A^{1/2}$: 
\begin{eqnarray}
P( | Z | > z ) & \le & 2  \exp \left (  - \frac{A}{B^2} h_1 \left ( \frac{Bz}{2A} \right ) \right )  \\
& = & 2  \exp \left (  - \frac{2}{L^2} h_1 \left ( \frac{Lz}{\tau} \right ) \right )\ \ \ \mbox{for all} \ \ z > 0.
\end{eqnarray}
\end{prop}
\smallskip

\par\noindent
{\bf Proof.}  
This follows from (\ref{hFcnInequalitiesTwo}) and elementary manipulations.
\hfill $\Box$
\medskip
 
 The classical route to proving inequalities of the form given in (\ref{BernsteinTypeInequality}) 
 for sums of independent random variables is via Bernstein's inequality; see for example 
 \cite{MR1385671} Lemmas 2.2.9 and 2.2.11, pages 102 and 103, or 
 \cite{MR3185193}, Theorem 2.10, page 37.   But the recent developments of concentration 
 inequalities via Stein's method yields inequalities of the form given in (\ref{BernsteinTypeInequality}) 
 for many random variables $Z$ which are {\sl not} sums of independent random variables: 
 see, for example, \cite{MR2763529, MR2773032} and \cite{MR3162712}.  The point of the previous 
 proposition is that (up to constants) these inequalities in terms of $h_0$ can be re-expressed in terms of the 
 (larger) function $h_1$. 
\medskip

\section{Bennett's inequality and the Bennett-Orlicz norm}
\label{sec:BennettOrliczNorm}

We begin with a statement of a version of Bennett's inequality for sums of
bounded random variables; 
see \cite{Bennett:62}, 
\cite{MR838963}, 
and 
\cite{MR3185193}.  
Let $h(x) \equiv x (\log x -1) +1$ and $h_2 (x) \equiv h(1+x)$.
This function arises in 
Bennett's inequality for bounded random variables and elsewhere; 
see e.g. 
\cite{Bennett:62}, 
\cite{MR838963}, 
and 
\cite{MR3185193}, 
page 35 (but note that their $h$ is our $h_2 = h(1+\cdot )$).
As noted in Example 1 above, the function $h$ also appears in exponential bounds for Poisson 
random variables:  
see 
\cite{MR838963}  
page 485, and 
\cite{MR3185193}  
page 23.
\medskip

\begin{prop}
\label{prop:Bennett-Ineq}
(Bennett)  
(i) \ Let $X_1, \ldots, , X_n$ be independent with 
$\max_{1 \le j \le n} (X_j - \mu_j ) \le b$, 
$E(X_j) = \mu_j$, $Var(X_j ) = \sigma_{n,j}^2$.  Let $\mu \equiv \sum_{j=1}^n \mu_j/n$, 
$\sigma_n^2 \equiv (\sigma_{n,1}^2 + \cdots + \sigma_{n,n}^2)/n$.  
Then with $\psi (x) \equiv 2 h(1+x)/x^2$,
\begin{eqnarray}
P \left (  \sqrt{n} ( \overline{X}_n - \mu ) \ge z \right )
& \le &  \exp \left ( - \frac{z^2}{2 \sigma_n^2} \psi \left ( \frac{z b}{\sqrt{n} \sigma_n^2  } \right ) \right ) \nonumber \\
& = & \exp \left ( - \frac{n \sigma_n^2}{b^2} h \left ( 1 + \frac{z b}{\sqrt{n} \sigma_n^2  } \right ) \right )  \nonumber \\
& = & \exp \left ( - \frac{n \sigma_n^2}{b^2} h_2 \left ( \frac{z b}{\sqrt{n} \sigma_n^2  } \right ) \right ) \label{BennettInequality}
\end{eqnarray}
for all $z > 0$.  \\
(ii) If,  in addition,  $\max_{1 \le j \le n} | X_j - \mu_j | \le b$, then
\begin{eqnarray*}
P \left (  | \sqrt{n} ( \overline{X}_n - \mu )| \ge z \right )
& \le &  2 \exp \left ( - \frac{z^2}{2 \sigma_n^2} \psi \left ( \frac{z b}{\sqrt{n} \sigma_n^2  } \right ) \right ) \\
& = & 2 \exp \left ( - \frac{n \sigma_n^2}{b^2} h \left ( 1 + \frac{z b}{\sqrt{n} \sigma_n^2  } \right ) \right ) \\
& = & 2 \exp \left ( - \frac{n \sigma_n^2}{b^2} h_2 \left ( \frac{z b}{\sqrt{n} \sigma_n^2 } \right ) \right ).
\end{eqnarray*}
\end{prop} 

Using the inequality $h(1+x) \ge 9 h_1 (x/3)$, it follows that 
\begin{eqnarray*}
 P \left ( | \sqrt{n} ( \overline{X}_n - \mu )| \ge z \right )
 \le 2\exp \left ( - \frac{9n \sigma_n^2}{b^2} h_1 \left ( \frac{z b}{3 \sqrt{n} \sigma_n^2  } \right ) \right ) .
 \end{eqnarray*}
Thus an inequality of the form of that in Lemma~\ref{lem:LvdG-lem1} holds 
with $2/L^2 = 9 n \sigma_n^2/b^2$ and $L/\tau = b/(3 \sqrt{n} \sigma_n^2 )$.  
Thus 
$L = \sqrt{2/9} b/(\sqrt{n} \sigma_n)$ and $\tau = L 3 \sqrt{n} \sigma_n^2 /b = \sqrt{2} \sigma_n$.  
It follows from Lemma~\ref{lem:LvdG-lem2} that 
$$
    \| \sqrt{n} ( \overline{X}_n - \mu ) \|_{\Psi_1 (\cdot ; \sqrt{3} L)} \le \sqrt{6} \sigma_n ,
$$
or
$$
\| \sqrt{n} ( \overline{X}_n - \mu ) \|_{\Psi_1 ( \cdot ; \sqrt{2/3} b / (\sqrt{n} \sigma_n ))} \le \sqrt{6} \sigma_n .
$$
But this bound has not taken advantage of the fact the the first bound above involves the 
function $h$ (or $h_2$) rather than $h_1$.   It would seem to be of potential interest 
to develop an Orlicz norm based on the function $ h_2 \equiv h(1+ \cdot ) $ rather than the function $h_1$. 
Motivated by the first inequality in 
Proposition~\ref{prop:Bennett-Ineq}, we define for each $L>0$ a new Orlicz norm based on the 
function $h_2$ as follows.
\begin{eqnarray*}
\Psi_2 (x; L) \equiv \exp \left ( \frac{2}{L^2} h_2(Lx) \right ) -1 .
\end{eqnarray*}
Since  $h_2$ is convex, $h_2 (0) = 0$, and $h_2$ is increasing on $[0,\infty)$, it follows that 
$\Psi_2 (\cdot; L)$ defines a valid Orlicz norm (as defined in Section~\ref{sec:intro}) for each $L$:
\begin{eqnarray}
\| X \|_{\Psi_2 (\cdot; L)} = \inf \left \{ c> 0 : \ E \Psi_2 \left ( \frac{|X|}{c} ; L\right ) \le 1  \right \} ;
\label{BennettOrliczNorm}
\end{eqnarray}
We call $\| X \|_{\Psi_2 (\cdot ;L)}$ the {\sl Bennett-Orlicz norm} of $X$.
Note that with $\psi (Lx) \equiv x^{-2} (2/L^2) h_2 (Lx)$,
\begin{eqnarray*}
\Psi_2 (x; L) & = & \exp \left ( x^2 \psi (Lx) \right ) -1 \\
& \sim & \left \{ \begin{array}{l l}  \exp (x^2) -1, & \mbox{for} \ L x \ \mbox{small},   \\
                                                    \exp \left ( \frac{2x}{L} \log (Lx) \right ) -1,   &  \mbox{for} \ L x \ \mbox{large} .
                        \end{array} \right . 
\end{eqnarray*}

We first relate $\Psi_2 (\cdot; L)$ to $\Psi_1 (x;L)$ and to the usual Gaussian Orlicz norm defined by 
$\Psi_2 (x) = \exp(x^2)-1$.

\begin{prop}
\label{prop:BoundsForBennettPsi} 
$\phantom{blab}$ \\
(i) \ $\Psi_2 (x; L) \le \exp (x^2) -1 = \Psi_2 (x)$ for all $x \ge 0$.\\
(ii) $\Psi_2 (x; L) \ge \Psi_1 (x ; L/3)  $ for $x \ge 0$.  
\end{prop}

\par\noindent
{\bf Proof.}  (i) follows since $\psi (x) \equiv 2 x^{-2} h(1+x) \le 1$ for all $x\ge 0$;  see 
\cite{MR838963}, 
Proposition 11.1.1, page 441.   To show that (ii) holds, note that by (\ref{PsiFcn-BernsteinOrlicz}) 
$$
\Psi_1 (x; L/3)  = \exp \left ( \frac{2 \cdot 9}{L^2} h_1 (Lx/3) \right ) -1 .
$$
Thus the claimed inequality in (ii) is equivalent to 
$$
\frac{2}{L^2 } h(1+Lx) \ge \frac{2\cdot 9}{L^2} h_1 ( Lx/3 ),
$$
or equivalently
$$
h(1+Lx) \ge 9 h_1 (Lx/3).
$$
But the inequality in the last display holds in view of (\ref{hFcnInequalitiesOne}).
\hfill $\Box$
\medskip

Note that while $h_1$  and $\Psi_1 (\cdot; L)$ have explicit inverses given in terms of $\sqrt{v}$ and $\log (1+v)$
by (\ref{InverseOFhOne}) and 
(\ref{BernsteinOrliczMaxIneq}), 
inverses of the functions $h_2$ and $\Psi_2 (\cdot ; L)$ can only be written in terms of Lambert's function 
(also called the {\sl product log} function) $W$ satisfying $W(z) \exp (W(z)) = z$; 
 see \cite{MR1414285}.  
 But this slight difficulty is easily overcome by way of several nice inequalities for $W$.
By use of $W$ and the 
inequalities developed in the Appendix,  Section 6, we obtain the following proposition concerning 
$\Psi_2^{-1} (\cdot ; L)$. 

\begin{prop} 
\label{BoundsForBennettInverse}   
$\phantom{blab}$ \\ 
(i)  \ $\Psi_2^{-1} (y; L) \le \Psi_1^{-1} (y; L/3) = \sqrt{\log (1+y)} + (L/6)\log (1+y) $ for $y \ge 0 $.\\
(ii) \ Furthermore, with $W$ denoting the Lambert $W$ function,
\begin{eqnarray*}
\Psi_2^{-1} (y; L) = \frac{1}{L}  h_2^{-1} \left ( \frac{L^2}{2} \log (1+y)  \right ) 
 =   \frac{1}{L} \left \{  \frac{\left ((L^2/2) \log (1+y) -1 \right )}{W\Big (\big ( (L^2/2) \log (1+y) -1 \big )/e \Big )} -1 \right \} .
\end{eqnarray*}
(iii) If $(L^2/2) \log (1+y) \ge 1$, then 
\begin{eqnarray*}
\Psi_2^{-1} (y; L ) & \le  &  L \frac{ \log (1+y)}{\log \left ( (L^2/2) \log (1+y) -1 \right ) } .
\end{eqnarray*}
(iv) If $(L^2/2) \log (1+y) \ge 5$, then 
\begin{eqnarray*}
\Psi_2^{-1} (y; L ) & \le & 2 L \frac{ \log (1+y)}{\log \left ( (L^2/2) \log (1+y) \right ) } \\
& \le & 2 L \frac{\log (1+y)}{\log \log (1+y)} \ \ \mbox{if also} \ \ L^2/2 \ge 1 .
\end{eqnarray*}
(v) If $(L^2/2)\log (1+y) \le 9/4$, then 
\begin{eqnarray*}
\Psi_2^{-1} (y; L ) & \le &  \sqrt{ 2\log (1+y)}.
\end{eqnarray*}
(vi) 
\begin{eqnarray*}
\Psi_2^{-1} (y; L ) 
& \le & \left \{ \begin{array}{l l}  (2.2/\sqrt{2}) \sqrt{\log (1+y)} & \ \ \mbox{if} \ (L^2/2) \log (1+y) \le 1+e, \\
                                        \ & \ \\ 
                                       L \frac{\log (1+y) -1}{\log ( (L^2/2) \log (1+y) -1 )} -1 & \ \ \mbox{if} \ \ (L^2/2)\log (1+y) > 1+e .
            \end{array} \right .
\end{eqnarray*}
\end{prop}

\medskip

\par\noindent
{\bf Proof.}  (i)  follows immediately from Proposition~\ref{prop:BoundsForBennettPsi}. 
(ii) follows from the definition of $\Psi_2 ( \cdot; L)$ and direct computation for the first part; 
the second part follows from Lemma~\ref{lem:h2InverseUsingLambertW}.  
The inequality in (iii)  follows from (ii) and Lemma~\ref{lem:LowerBoundForW}.    
The first inequality in (iv) follows from (iii) since $\log (y-1) \ge (1/2) \log y$ for $y \ge 4$.  
The second inequality in (iv) follows by noting that 
$$
\log ((L^2/2)\log (1+y)) = \log (L^2/2) + \log \log (1+y) \ge \log\log (1+y)
$$
if $L^2/2 \ge 1$.  
(v) follows from (ii) and Lemma~\ref{lem:UpperBoundsForHInverseandH2Inverse}, part (iv). 
\hfill $\Box$
\medskip

Lemmas~\ref{lem:LvdG-lem1} and~\ref{lem:LvdG-lem2}  by \cite{MR3101846} 
as stated in Section~\ref{BernsteinOrlicz} should 
be compared with the development on
page 96 of \cite{MR1385671}.  
We now show that the following analogues of Lemmas~\ref{lem:LvdG-lem1} - ~\ref{lem:LvdG-lem3}  hold for $\| Z \|_{\Psi_2 (\cdot; L)}$.
\medskip
 
\begin{lem}
\label{lem:1A}
Let $\tau \equiv \| Z \|_{\Psi_2 (\cdot; L)}$.  Then 
$$
P\left ( | Z | > \frac{\tau}{L} h_2^{-1} (L^2 t/2) \right ) \le 2 e^{-t}  \ \ \mbox{for all} \ \ t>0 
$$
where $h_2 (x) \equiv h(1+x)$ and $h_2^{-1}$ is the inverse of $h_2$ (so that 
$h_2^{-1} (y) = h^{-1} (y) -1$).
\end{lem}
\smallskip

\par\noindent
{\bf Proof.}  Let $y>0$.  Since $\Psi_2 (x;L) = \exp ((2/L^2) h_2(Lx))-1 = e^t-1$ implies $h_2(Lx) = L^2t/2$, 
it follows that for any $c> \| Z \|_{\Psi_2 (\cdot; L)}$ we have 
\begin{eqnarray*}
P\left ( \frac{|Z|}{c} > \frac{1}{L} h_2^{-1} (L^2t/2) \right ) 
& = & P \left (  h_2 \left ( \frac{L|Z|}{c}\right ) > L^2 t/2  \right ) \\
& = & P\left ( \frac{2}{L^2} h_2 \left ( \frac{L |Z|}{c}\right ) >  t \right ) \\
& = & P\left ( \exp \left ( \frac{2}{L^2} h_2 \left ( \frac{L |Z|}{c}\right ) \right ) -1 > e^t-1 \right ) \\
& = & P \left ( \Psi_2 \left (  \frac{|Z|}{c}; L \right ) > e^t -1 \right ) \\
& \le & \frac{E \left \{ \Psi_2 \left ( \frac{|Z|}{c}; L \right ) +1 \right \} }{e^t}  \\
& \rightarrow & 2 e^{-t} \ \ \mbox{as} \ c \downarrow \tau \equiv \| Z \|_{\Psi_2 (\cdot; L)} .
\qquad \qquad \qquad \qquad \qquad \Box
\end{eqnarray*}
 
\begin{lem}
\label{lem:2A}
Suppose that for some $\tau >0$ we have
\begin{eqnarray*}
P\left ( | Z | > \frac{\tau}{L} h_2^{-1} (L^2 t/2) \right ) \le 2 e^{-t}  \ \ \mbox{for all} \ \ t>0.
\end{eqnarray*}
Equivalently,
\begin{eqnarray*}
P( | Z | > z ) 
& \le & 2 \exp \left ( - \frac{2}{L^2} h_2  \left ( \frac{Lz}{\tau} \right ) \right ) \\
& = & 2 \exp \left ( - \frac{2}{L^2} h  \left (1+ \frac{Lz}{\tau} \right ) \right )\\
& = & 2 \exp \left ( - \frac{z^2}{\tau^2} \psi \left ( \frac{Lz}{\tau} \right ) \right ) \ \ \mbox{for all} \ \ z >0 .
\end{eqnarray*}
Then $\| Z \|_{\Psi_2(\cdot ; \sqrt{3}L )} \le \sqrt{3} \tau $.
\end{lem}
\medskip

\par\noindent
{\bf Proof.} Let $\alpha, \beta >0$.   We compute
\begin{eqnarray*}
E \Psi_2 \left ( \frac{|Z|}{\alpha \tau }; \beta L \right ) 
& = & \int_0^{\infty} P\left ( \Psi_2 \left ( \frac{|Z|}{\alpha \tau}; \beta L \right )  \ge v \right ) dv  \\
& = & \int_0^{\infty} P\left ( \frac{2}{\beta^2 L^2} h_2 \left (\frac{\beta L |Z|}{\alpha \tau} \right ) \ge \log (1+v) \right ) dv \\
& = & \int_0^{\infty} P \left ( \frac{ \beta L| Z |}{\alpha \tau}  \ge \tau h_2^{-1} \left ( \frac{\beta^2 L^2}{2} \log (1+v) \right ) \right ) dv \\
& = & \int_0^{\infty} P\left ( | Z | \ge \frac{\alpha \tau}{\beta L} h_2^{-1} \left ( \frac{\beta^2 L^2}{2} \log (1+v) \right ) \right )dv \\
& = & \int_0^{\infty} P \left ( | Z | \ge \frac{\alpha \tau}{\beta L} h_2^{-1} \left ( \frac{\beta^2 L^2}{2} t \right ) \right ) e^{t} dt .
\end{eqnarray*}
Choosing $\alpha = \beta = \sqrt{3}$ this yields
\begin{eqnarray*}
E \Psi_2 \left ( \frac{|Z|}{\sqrt{3} \tau }; \sqrt{3} L \right ) 
& = &  \int_0^{\infty} P \left ( | Z | \ge \frac{ \tau}{ L} h_2^{-1} \left ( \frac{ L^2}{2} 3 t \right ) \right ) e^{t} dt \\
&\le &  \int_0^{\infty} 2 \exp (-3 t) \exp (t) dt  = 1 .
\end{eqnarray*}
Hence we conclude that 
$\| Z \|_{\Psi_2 (\cdot ; \sqrt{3} L)} \le \sqrt{3} \tau $. 
\hfill $\Box$
\medskip 

\begin{cor} 
\label{cor:BennettOrliczNormsPoisson-Binom-Gauss}
$\phantom{blab}$\\
(i) \ If $X \sim \mbox{Poisson} (\nu)$, then $\| X - \nu \|_{\Psi_2 (\cdot; \sqrt{6/\nu} )} \le \sqrt{6 \nu }$. \\
(ii) \ If $X_1, \ldots , X_n $ are i.i.d. Bernoulli$(p)$, then
$$
 \| \sqrt{n} (\overline{X}_n  - p )\|_{\Psi_2 (\cdot; \sqrt{2}/(\sqrt{n p(1-p)} )) } \le  \sqrt{6 p(1-p)} .
$$
(iii) \ 
If $X \sim N(0,1)$, then $\| X \|_{\Psi_h (\cdot ; L)} \le \sqrt{6} $ for every $L>0$.
By taking the limit on $L \searrow 0$ 
and noting that $\Psi_2 (z ; L ) \rightarrow  \Psi_2 (z )\equiv \exp (z^2) -1$ as $L \searrow 0$
 this yields $\|  X \|_{\Psi_2 } \le \sqrt{6} $.
In this case it is known that $\| X \|_{\Psi_2} = \sqrt{8/3}$.  (See \cite{MR1385671}, Exercise 2.2.1, page 105.)
\end{cor}
\medskip

Now for an inequality paralleling Lemma~\ref{lem:LvdG-lem3} for the Bernstein-Orlicz norm:
\smallskip
 
\begin{lem}
\label{lem:3A}
Let $\tau $ and $L$ be constants, and let $Z_1, \ldots , Z_m$ 
be random variables satisfying $\max_{1 \le j \le m} \| Z_j \|_{\Psi_2( \cdot ; L)} \le \tau$. 
Then
\begin{eqnarray*}
E \{ \max_{1 \le j \le m} | Z_j | \} 
& \le & \tau \Psi_2^{-1} (m ; L )\\
& \le & 2 \tau L \frac{\log (1+m)}{ \log \log (1+m) }  \ \ \mbox{if} \ \ L^2 \ge 2 \ \mbox{and} \ \log (1+m) \ge 5 .
\end{eqnarray*}
Furthermore,
\begin{eqnarray*}
E \{ \max_{1 \le j \le m} | Z_j | \} 
& \le & \tau \Psi_2^{-1} (m ; L )
             \le 2 \tau  \left \{ L \frac{\log (1+ m)}{\log \log (1+m)}  + \sqrt{ \log (1+ m)}  \right \} 
\end{eqnarray*}
for all $m$ such that $\log(1+m) \ge 5$ (or $m \ge e^5-1$).
\end{lem}
\smallskip
 
\begin{rem}
The point of this last bound is that it gives an explicit trade-off between 
the Gaussian component (the term $\sqrt{\log(1+m)}$) and the Poisson component (the term
$\log (1+m) / \log \log (1+m)$) governed by a Bennett type inequality.  
In contrast, the bounds obtained by 
\cite{MR3101846} yield a trade-off between the Gaussian world and the sub-exponential world
governed by a Bernstein type inequality.
\end{rem}
\medskip

\par\noindent
{\bf Proof.}  We write $\Psi_{2,L} \equiv \Psi_2 (\cdot ; L)$.  Let $c > \tau$.  Then by Jensen's inequality
\begin{eqnarray*}
E \left \{ \max_{1 \le j \le m} | Z_j | \right \} 
& \le & c \Psi_{2,L}^{-1} \left ( E \{ \Psi_{2,L} ( \max_{1 \le j \le m} | Z_j | /c ) \} \right )\\
& = & c \Psi_{2,L}^{-1} \left ( E \{ \max_{1 \le j \le m} \Psi_{2,L} ( | Z_j |/c ) \} \right ) \\
& \le & c \Psi_{2,L}^{-1} \left ( \sum_{j=1}^m E \Psi_{2,L}  ( | Z_j | /c ) \right ) \\
& \le & c \Psi_{2,L}^{-1} \left ( m \max_{1 \le j \le m} E \Psi_{2,L}  ( | Z_j | /c ) \right ) .
\end{eqnarray*}
Therefore, 
\begin{eqnarray}
E \left \{  \max_{1 \le j \le m} | Z_j | \right \}  
& \le & \lim_{c \searrow \tau} \left \{ c \Psi_{2,L}^{-1} \left ( m \max_{1 \le j \le m} E \Psi_{2,L}  ( | Z_j | /c ) \right ) \right \} \\
& = & \tau \Psi_{2,L}^{-1} ( m ) = \frac{\tau}{L} h_2^{-1} \left ( \frac{L^2}{2} \log (1+m) \right ) .
\label{MaxBoundInTermsOfhInverse}
\end{eqnarray} 
The remaining claims follow from Proposition~\ref{BoundsForBennettInverse}.  
\hfill $\Box$
\bigskip

Here are analogues of Lemmas 4 and 5 of \cite{MR3101846}.
\medskip

\begin{lem}
\label{lem:Lemma4A}
Let $Z_1, \ldots , Z_m$ be random variables satisfying 
\begin{eqnarray}
\max_{1 \le j \le m} \| Z_j \|_{\Psi_2 (\cdot, L)} \le \tau
\label{MaxOfOrliczNormsBound}
\end{eqnarray}
for some $L$ and $\tau$.  Then, for all $t>0$
\begin{eqnarray*}
P \left ( \max_{1 \le j \le m} | Z_j | \ge \frac{\tau}{L} \left ( 
  h_2^{-1} ( L^2 t/2) + h_2^{-1} (L^2 \log (1+m)/2) \right ) \right ) \le 2 e^{-t} .
\end{eqnarray*}
\end{lem}
\medskip

\par\noindent
{\bf Proof.}    For any $a>0$ and $t>0$ concavity of $h_2^{-1}$ together with $h_2^{-1} (0) = 0$ imply that 
\begin{eqnarray*}
h_2^{-1} (a) + h_2^{-1} (t) \ge h_2^{-1} (a+t) .
\end{eqnarray*}
Therefore, by using a union bound and Lemma~\ref{lem:1A}
\begin{eqnarray*}
\lefteqn{P \left ( \max_{1 \le j \le m} | Z_j | \ge \frac{\tau}{L} \left ( 
  h_2^{-1} ( L^2 t/2) + h_2^{-1} (L^2 \log (1+m)/2) \right ) \right )} \\
& \le & P \left ( \max_{1 \le j \le m} | Z_j | \ge \frac{\tau}{L} \left ( 
  h_2^{-1} ( 2^{-1} L^2 ( t +  \log (1+m))) \right ) \right )\\
& \le & \sum_{j=1}^m P\left ( | Z_j | \ge \ \frac{\tau}{L} \left ( 
          h_2^{-1} ( 2^{-1} L^2 ( t +  \log (1+m))) \right ) \right ) \\
& \le & 2m \exp \left ( - (t + \log (1+m))\right ) = 2 \frac{m}{m+1} e^{-t}  \le 2 e^{-t}.
\end{eqnarray*}
\hfill $\Box$
\medskip

\begin{lem}
\label{lem:Lemma5A}
Let $Z_1, \ldots , Z_m$ be random variables satisfying 
(\ref{MaxOfOrliczNormsBound}).  Then
\begin{eqnarray*}
\bigg \| \left ( \max_{1 \le j \le m} | Z_j | - \frac{\tau}{L} h_2^{-1} \left ( \frac{L^2}{2} \log (1+m) \right ) \right )_{+} 
\bigg \|_{\Psi_2 (\cdot ; L)} \le \sqrt{3} \tau .
\end{eqnarray*}
\end{lem}

\par\noindent
{\bf Proof.}  Let 
\begin{eqnarray*}
Z \equiv \left ( \max_{1 \le j \le m} | Z_j | - \frac{\tau}{L} h_2^{-1} \left ( \frac{L^2}{2} \log (1+m) \right ) \right )_{+}  .
\end{eqnarray*}
Then Lemma~\ref{lem:Lemma4A} implies that 
\begin{eqnarray*}
\lefteqn{P\left ( Z \ge \frac{\tau}{L} h_2^{-1} \left ( \frac{L^2 t}{2} \right ) \right )} \\
& \le & P \left ( \max_{1 \le j \le m} | Z_j | \ge \frac{\tau}{L} \left ( h_2^{-1} \left ( \frac{L^2}{2} \log (1+m) \right ) 
                 + h_2^{-1} \left ( \frac{L^2 t}{2} \right ) \right ) \right ) \\
& \le & 2 e^{-t} .
\end{eqnarray*}
Then the conclusion follows from Lemma~\ref{lem:2A}.  
\hfill $\Box$

\section{Prokhorov's ``arcsinh'' exponential bound and Orlicz norms} 
\label{sec:ProkhorovInequality}

Another important exponential bound for sums of independent bounded random variables is 
due to \cite{MR0121857}.  As will be seen below, Prokhorov's bound involves another function $h_4$ 
(rather than $h_2$ of Bennett's inequality) given by 
\begin{eqnarray}
h_4 (x)  = (x/2) \mbox{arcsinh}(x/2) = (x/2) \log \left (x/2 + \sqrt{1+ (x/2)^2} \right ) .
\label{h4ProhorovDefn}
\end{eqnarray}
Suppose that $X_1, \ldots , X_n$ are independent random variables with $E(X_j) = \mu_j$ and 
$|X_j -\mu_j | \le b$ for some $b>0$.  Let $S_n = X_1 + \cdots + X_n$, 
and set $\mu \equiv n^{-1} \sum_{j=1}^n \mu_j$, $\sigma_n ^2 \equiv n^{-1} Var(S_n)$.
Prokhorov's ``arcsinh'' exponential bound is as follows:  
\medskip

\begin{prop}
\label{prop:ProkhorovIneq}
(Prokhorov) If the $X_j$'s satisfy the above assumptions, then 
\end{prop}
\vspace{-.1in}
\begin{eqnarray*}
P( S_n - n \mu \ge z ) \le \exp \left ( - \frac{z}{2 b} \mbox{arcsinh} \left ( \frac{z b}{2 \sigma_n^2} \right ) \right ) .
\end{eqnarray*}
{\it Equivalently, with $\sigma_n^2 \equiv n^{-1} Var(S_n)$ and} $h_4 (x) \equiv (x/2) \mbox{arcsinh} (x/2)$,
 \begin{eqnarray}
P( |\sqrt{n} (\overline{X}_n - \mu) |  \ge z ) 
& = & 2\exp \left ( - \frac{z \sqrt{n} }{2 b}  \mbox{arcsinh} \left ( \frac{z  b}{2 \sqrt{n} \sigma_n^2} \right ) \right ) \nonumber \\
& = & 2\exp \left ( - \frac{ n \sigma_n^2}{b^2}  \left (\frac{z b}{2 \sqrt{n} \sigma_n^2} \right ) 
                              \mbox{arcsinh} \left ( \frac{z  b}{2 \sqrt{n} \sigma_n^2} \right ) \right ) \nonumber \\
& \equiv & 2 \exp \left ( - \frac{ n \sigma_n^2}{b^2} h_4 \left ( \frac{z b}{\sqrt{n} \sigma_n^2} \right ) \right ) .
\label{ProhorovArcSinhInequality}
\end{eqnarray}
\medskip

See e.g.
\cite{MR0121857},   
\cite{MR0455094},   
\cite{MR1666908},   
\cite{MR770640},  
and \cite{MR2331988t} . 
\cite{MR770640} use Prokhorov's inequality to control Orlicz norms for  functions $\Psi$ of 
the form $\Psi (x) = \exp ( \psi (x))$ with $\psi (x) \equiv x\log(1+x)$ and use the resulting inequalities
to show that the optimal constants $D_p$ in Rosenthal's inequalities grow as $p/\mbox{log}(p)$.

\medskip

\par\noindent
\cite{MR2331988t} gives an improvement of Prokhorov's inequality which involves replacing 
$h_4 $ by 
$$
h_5(x) \equiv x \mbox{arcsinh} (x/2) - 2 \Big ( \mbox{cosh} \big (\mbox{arcsinh}(x/2) \big ) -1 \Big ) .
$$
Note that Prokhorov's inequality is of the same form as Bennett's inequality (\ref{BennettInequality}) 
in Proposition~\ref{prop:Bennett-Ineq}, 
but with Bennett's $h_2$ replaced by Prokhorov's $h_4$.
\medskip

\par\noindent
Thus we want to
 compare Prokhorov's inequality (and Kruglov's improvement thereof) to Bennett's inequality.
As can be seen from the above development, this boils down to comparison of 
the functions $h_2$, $h_4$, and $h_5$.    
The following lemma makes a number of comparisons and contrasts between the functions $h_2$, $h_4$, and $h_5$.\\

\begin{lem} 
\label{Prohorov-Bennett-Comparison}
(Comparison of $h_2$, $h_4$, and $h_5$)\\
(i)(a) \ \ \ \ \ \ $h_2(x) \ge h_5 (x) \ge h_4 (x) $ for all $x \ge 0$. \\
(i)(b) \ \ \ \ \ $h_2^{-1} (y) \le h_5^{-1} (y) \le h_4^{-1} (y) $ for all $y \ge 0$. \\ 
$\phantom{blab}$\\
(ii)(a) \ \ $h_2 (x) \ge (x/2) \log (1+x) \ge (x/2) \log (1 + x/2)$ for all $x \ge 0$.\\
(ii)(b) \ \ $h_4 (x) \ge (x/2) \log (1 + x/2)$ for all $x \ge 0$.\\
(ii)(c) \ \ $h_5 (x) \ge \ (x/2)\log (1+x/2) $ for all $x\ge 0$.\\
$\phantom{blab}$\\
(iii)(a) \ $h_2 (x) \sim 2^{-1} x^2$ as $x \searrow 0$;   $h_2 (x) \sim  x \log (x)$ as $x\rightarrow \infty$.\\
(iii)(b) \ $h_4(x) \sim 4^{-1} x^2 $ as $x \searrow 0$;   $h_4 (x) \sim (1/2) x \log (x)$  as $x\rightarrow \infty$.\\
(iii)(c) \ $h_5 (x) \sim 4^{-1} x^2$ as $x \searrow 0$;   $h_5 (x) \sim  x \log (x)$  as $x\rightarrow \infty$.\\
(iii)(d) \ $h_2(x) - h_4 (x) \sim x^2/4$ as $x \searrow 0$;   $h_2 (x) - h_4 (x) \sim (1/2) x\log x$ as $x \rightarrow \infty$. \\
(iii)(e) \ $h_2(x) - h_5 (x) \sim x^2/4$ as $x \searrow 0$;   $h_2(x) - h_5 (x) \sim \log x  $ as $x \rightarrow \infty $. 

$\phantom{blab}$\\
(iv)(a) \ $h_2 (x) = 2^{-1} x^2 \psi_2 (x) $ where 
$$
\psi_2 (x)  \ge \left \{ \begin{array}{l} x^{-1} \log (1+x) , \\ (1+x/3)^{-1} .\end{array} \right .
$$
(iv)(b) \ $h_4 (x) = 4^{-1} x^2 \psi_4 (x)$ where 
$$
\psi_4 (x)  \ge \left \{ \begin{array}{l l} 2x^{-1} \log (1+x/2) , & \mbox{for} \ x \ge 0  \\ 
                                                          (1/2)/(1+x/2)^{-1} , & \mbox{for} \ x \ge 0 \\ 
                                                          (1-\delta)/(1+x/2), &  \mbox{for}\ x \le 2\delta^{1/2}/(1/2-\delta)^{1/2} .
       \end{array}  \right .
$$
(iv)(c) \ $h_5 (x) = 4^{-1} x^2 \psi_5 (x)$ where 
$$
\psi_5 (x) \ge \left \{ \begin{array}{l l} 2x^{-1} \log (1+x/2) , &  \mbox{for} \ x \ge 0  \\ 
                                                          1/(1+x/2)^{-1} , & \mbox{for} \ x \ge 0 .
                                 \end{array}  \right .
$$
\end{lem}
\medskip

\par\noindent
{\bf Proof.}  (i)  We first prove that $h_2 (x) \ge h_4 (x)$.  Let $g(x) = h_2 (x) - h_4 (x)$; thus
\begin{eqnarray*}
g(x) = (1+x)\log (1+x) - x - (x/2) \log (x/2 + \sqrt{1 + (x/2)^2}) .
\end{eqnarray*}
Then $g(0) = 0$ and
\begin{eqnarray*}
g'(x) = \log (1+x) - \frac{1}{2} \log \left ( x/2 + \sqrt{1 + (x/2)^2}  \right ) - \frac{x/4}{\sqrt{1 + (x/2)^2}} 
\end{eqnarray*}
also has $g' (0) = 0$.  Note that $\sqrt{1+(x/2)^2} \le 1 + x/2$ and hence $x/2 + \sqrt{1+(x/2)^2} \le 1 + x$.
Thus
\begin{eqnarray}
\log \big ( (x/2) + \sqrt{1 + (x/2)^2} \big ) \le \log (1+x),
\label{ArcSinh-Log-Inequality}
\end{eqnarray}
and hence
\begin{eqnarray*}
g'(x) 
& \ge & \log (1+x) - (1/2) \log (1+x) -  \frac{x/4}{\sqrt{1 + (x/2)^2}} \\
& = & (1/2) \log (1+x) - \frac{x/4}{\sqrt{1 + (x/2)^2}},
\end{eqnarray*}
and it suffices to show that the right side is $\ge 0$ for all $x$. Thus we let 
$$
m(x) \equiv g'(x) = \frac{1}{2} \log (1+x) - \frac{1}{2} \frac{x/2}{\sqrt{1+ (x/2)^2}} 
= \frac{1}{2} \left \{ \log (1+x) - \frac{x/2}{\sqrt{1+(x/2)^2}} \right \} .
$$
Let $\overline{m} (x) \equiv 2 m (2x) = \log (1+2x) - \frac{x}{\sqrt{1+x^2}} $.
Then $\overline{m} (0) = 0$ and we compute 
\begin{eqnarray*}
\overline{m}' (x) 
& = & \frac{2}{1+2x} - \frac{1}{\sqrt{1+x^2}}  +  \frac{x^2}{(1 + x^2)^{3/2}}  \\
& = & \frac{2}{1+2x} - \frac{1}{\sqrt{1+x^2}} \left (1- \frac{x^2}{(1+x^2)} \right )\\
& = & \frac{2 (1+x^2)^{3/2} - (1+2x)}{(1+2x)(1+x^2)^{3/2}} \equiv \frac{j(x)}{(1+2x)(1+x^2)^{3/2}} 
\end{eqnarray*}
so that $\overline{m}' (0) = 1$ and the numerator, $j$, is easily seen to be non-negative 
since $(1+x^2)^{3/2} \ge 1+x^2$ implies $2 (1+x^2)^{3/2} \ge 2 (1+x^2) \ge 1+2x$  for all $x \ge 0$. 
Thus $h_2 (x) \ge h_4 (x) $.  

\cite{MR2331988t} shows that $h_5 (x) \ge h_4 (x)$.  Now we show that $h_2 (x) \ge h_5 (x)$. 
Note that with $g(x) \equiv h_2 (x) - h_5 (x)$, 
$$
g' (x) = \log (1+x) - \mbox{arcsinh}(x/2) 
$$ 
has $g'(x) = 0$ and $g'(x) \ge 0$ (as was shown above in (\ref{ArcSinh-Log-Inequality}) ). 
Thus $g(x) = \int_0^x g'(v) dv \ge 0$. 

(i)(b)  The inequalities for the inverse functions follow immediately from the inequalities for the functions 
themselves in (i)(a).    
\medskip

\par\noindent
(ii)(a) To show that the first inequality holds, consider 
\begin{eqnarray*}
g(x) & \equiv & h_2 (x) - (x/2) \log (1+x) \\
& = & (1+x) \log (1+x) -x - (x/2) \log (1+x) \\
& = & (1+x/2) \log (1+x) - x .
\end{eqnarray*}
Then $g(0) = 0$ and 
\begin{eqnarray*}
g' (x) & = & \frac{1}{2} \log (1+x) + \frac{1+x/2}{1+x} -1 \\
& = & \frac{1/2}{1+x} \left \{ (1+x) \log (1+x) - x \right \} = \frac{1/2}{1+x} h_2 (x) \ge 0 .
\end{eqnarray*}
Thus $g' (0) = 0$ and $g(x) = \int_0^x g' (y) dy \ge 0$.  The second inequality in (ii)(a) is trivial.\\

\par\noindent
(ii)(b) \ This follows easily from $\mbox{arcsinh}(v) = \log(v+\sqrt{1+v^2}) \ge \log (v+1)$ for all $v \ge 0$.\\
(ii)(c) \ This follows from (i)(a) and (ii)(b).
\smallskip

\par\noindent
(iii)(a) This follows from $\psi_2 (x) \equiv \psi (x) \rightarrow 1$ as $x\searrow 0$; see Proposition 11.1.1, page 441, 
\cite{MR838963}.\\
(iii)(b) Now 
\begin{eqnarray*}
h_4^{\prime} (x) = \frac{x}{4 \sqrt{1 + (x/2)^2}} + \frac{1}{2} \mbox{arcsinh} (x/2) ,  
\end{eqnarray*}
with $h_4 ' (0) = 0$, and 
\begin{eqnarray*}
h_4^{\prime \prime} (x) = \frac{1}{2 \sqrt{1 + (x/2)^2}} - \frac{x^2}{16 (1 + (x/2)^2)^{3/2}} = \frac{8 + x^2}{2 (4 + x^2)^{3/2}} 
\end{eqnarray*}
with $h_4^{\prime \prime} (0) = 1/2$.  Therefore 
$$
h_4 (x) = h_4^{\prime \prime} (x^*) \frac{x^2}{2} \ \ \mbox{for some} \ \ 0 \le x^* \le x
$$
and 
$$
\frac{4}{x^2} h_4 (x) = 2 h_4^{\prime \prime} (x^*) \rightarrow 1 \ \ \mbox{as} \ \ x \ge x^* \searrow 0.
$$ 
(iii)(c) \ Now 
\begin{eqnarray*}
h_5^{\prime} (x) & = & \mbox{arcsinh}(x/2) , \\
h_5^{\prime \prime} (x) & = & \frac{1}{2 \sqrt{1 + (x/2)^2}} \rightarrow \frac{1}{2} \ \ \mbox{as} \ x \searrow 0,
\end{eqnarray*}
where $h_5^{\prime} (x) = 0$ and $h_5^{\prime \prime}$ is decreasing.
Thus $h_5 (x) = (x^2/2) h_5^{\prime \prime} (x^*) $ for some $0 \le x \le x^*$ and we conclude that 
$4 x^{-2} h_5 (x) \rightarrow 1$ as $x \ge x^* \searrow 0$. 
\medskip

\par\noindent 
(iv)(a) \ The first part is a restatement of (ii)(a).   The second part follows from (\ref{hFcnInequalitiesTwo}):  
$h_2(x) = h(1+x) \ge 9 h_0 (x) = x^2/(2(1+x/3))$, and the claim follows by definition of $\psi_2$. \\
(iv)(b) The first inequality is a restatement of (ii)(b).  The second inequality follows since 
$h_4 (x) = h_4^{\prime \prime} (x^*) $ where $x \mapsto h_4^{\prime \prime} (x)$ is decreasing,
so 
\begin{eqnarray*}
\frac{4}{y^2} h_4 (x) 
& = & 2 h_4^{\prime \prime} (x^*) \ge 2 h_4^{\prime \prime} (x) \\
& = & \frac{1}{\sqrt{1 + (x/2)^2}} - \frac{x^2}{8 (1 + (x/2)^2)^{3/2}} 
          = \frac{1}{\sqrt{1+ (x/2)^2}} \cdot \frac{1+x^2/8}{1+x^2/4} \\
 & \ge & \frac{1/2}{\sqrt{1+(x/2)^2}} \ge \frac{1/2}{1 + x/2} .
\end{eqnarray*}
To prove the third inequality, note that 
\begin{eqnarray*}
\frac{1+x^2/8}{1+x^2/4} \ge c 
\end{eqnarray*}
holds if $1 + x^2/8 \ge c (1+ x^2/4) $, or if $1-c \ge (x^2/4) (c-1/2)$.  
Then rearrange and take $c = (1-\delta)$ for $\delta \in (0,1/2)$.\\
(iv)(c) The first inequality follows from (ii)(c).  The second inequality follows by
arguing as in (iv)(b),  but now without the complicating second factor: note that
\begin{eqnarray*}
\frac{4}{x^2} h_5 (x)  = 2 h_5^{\prime \prime} (x^*) 
 \ge  2 h_5^{\prime \prime} (x) = \frac{1}{\sqrt{1 + (x/2)^2}}
 \ge   \frac{1}{1 + x/2} 
\end{eqnarray*}
since $h_5^{\prime\prime}$ is decreasing.  
\hfill  $\Box$
\medskip 

\par
\begin{figure}[ht]
    \centering
    \includegraphics[width=\linewidth,height=5.5cm,keepaspectratio]{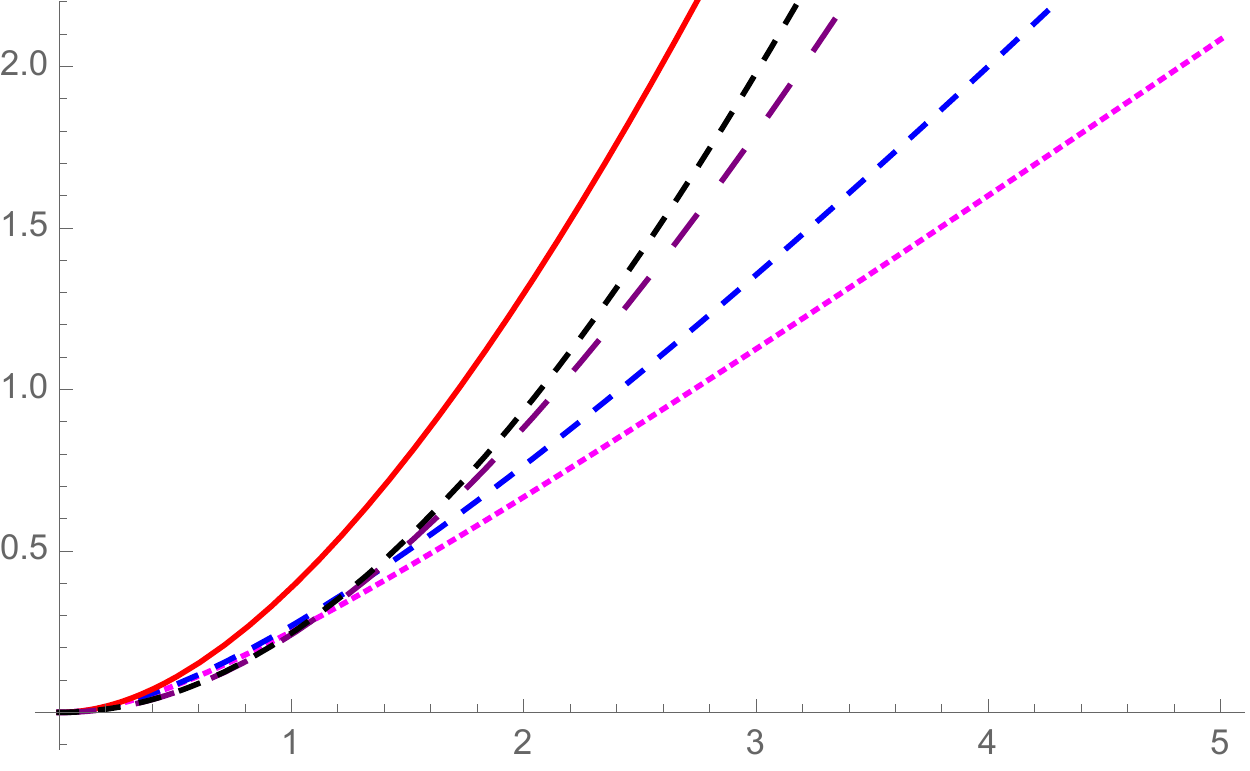}
    \caption{Comparison of the $h_k$ functions $h_0$, $h_1$, $h_2$, $h_4$, and $h_5$.
                  The  plot shows the functions $h_k$. 
                  The function $h_0$ is plotted in magenta (tiny dashing), 
                  $h_1 $ in blue (medium dashing),  
                  $h_2 $ in red (no dashing), 
                  $h_4$ in purple (large dashing), 
                  and $h_5$ is plotted in black (medium dashing).  
                  For values of the argument larger than $\approx 1.4$ $h_2> h_5 > h_4 >> h_1> h_0$ (and all are below $h_2$),
                  while for values of the argument smaller than $\approx 1$,  $h_2 > h_1 > h_0 >> h_5 > h_4$.}
     \label{fig:fig1}
 \end{figure}

\begin{figure}[ht]
    \centering
    \includegraphics[width=\linewidth,height=5.5cm,keepaspectratio]{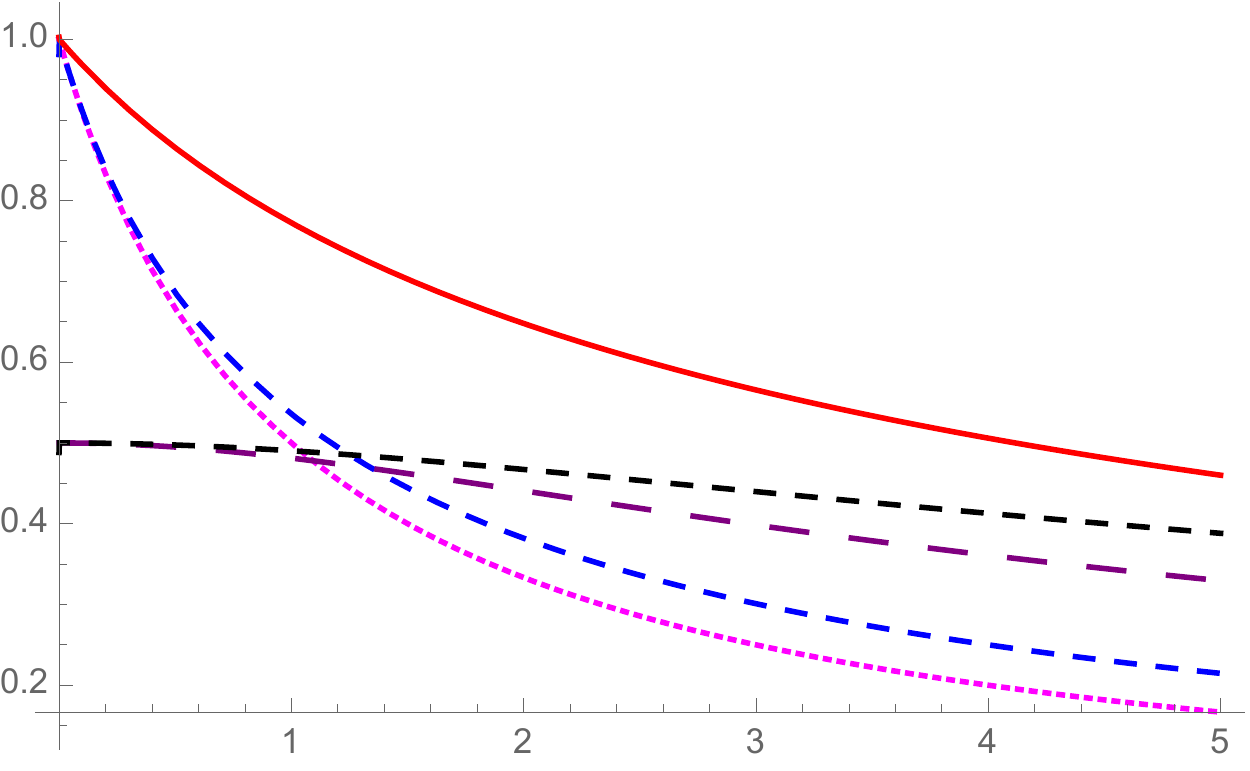} 
    \caption{The plot depicts (with the same colors and dashing as in Figure~\ref{fig:fig1}) the ratios 
             $x \mapsto h_k(x) /(x^2/2) \equiv \psi_k (x) $ for $k \in \{ 0,1,2,4,5 \}$. 
             This figure illustrates our finding that the Prokhorov type $h-$functions are smaller by a factor of $1/2$ at $x=0$,
              while they again dominate the Bernstein type $h-$functions for larger values of $x$, with the cross-overs 
              occurring again between $1$ and $1.4$.}
    \label{fig:fig2}
 \end{figure}

\noindent
{\bf Discussion:}  
{\bf 1.} \ Even though Kruglov's inequality improves on Prokhorov's inequality, 
(ia) of Lemma~\ref{Prohorov-Bennett-Comparison} 
shows  that Bennett's inequality dominates both  Kruglov's improvement 
of Prokhorov's inequality and Prokhorov's inequality itself:  $h_2 \ge h_5 \ge h_4$. \\
{\bf 2.}  \ (ii) of Lemma~\ref{Prohorov-Bennett-Comparison} shows that all three of the inequalities, Bennett, Kruglov, and Prokhorov,
are based on functions $h_2$, $h_5$, and $h_4$ which are bounded below by $(x/2)\log (1+x/2)$ for all $x \ge 0$. 
On the other hand,
(ii)(d) shows that both $h_2$ and $h_5$ are very nearly equivalent for large $x$, but that 
although $h_4 $ grows at the same $x\log x$ rate as $h_2 $ and $h_5$, $h_4$ is smaller by a multiplicative factor of $1/2$ 
as $x \rightarrow \infty$. \\
{\bf 3.} \  (iii)(a-c) of Lemma~\ref{Prohorov-Bennett-Comparison} 
shows that $h_2 (x) \sim x^2/2$ as $x \searrow 0$ while $h_k (x) \sim x^2/4$ for both $h_5 $ and $h_4$; 
thus $h_2(x)$ is larger at $x=0$ by a factor of $2$. 
 Furthermore, the difference  $h_2 - h_4$  is of order $(1/2)x \log x$ as $x\rightarrow \infty$, while the difference
 $h_2 - h_5$ is only of order $\log x$ as $x\rightarrow \infty$.  \\
 {\bf 4.} \ (iv) of Lemma~\ref{Prohorov-Bennett-Comparison} re-expresses the behavior 
 of the Kruglov and Prokhorov inequalities for small values of $x$ 
 in terms of the corresponding $\psi_k$ functions.\\
 The upshot of all of these comparisons is that Bennett's inequality dominates both the Kruglov and Prokhorov inequalities. 
 Figures~\ref{fig:fig1}~-~\ref{fig:fig2} give graphical versions of these comparisons as well as comparisons to the 
 Bernstein type $h-$functions $h_0$ and $h_1$.
\medskip 

\section{Comparisons with some results of Talagrand}
\label{sec:TalagrandComparisons}

Our goal in this section is to give comparisons with some results of 
\cite{MR1048946} and 
\cite{MR1258865}, 
especially his Theorem 3.5, page 45, and Proposition 6.5, page 58.  

\cite{MR1258865} 
defines a function $\varphi_{L,S}$ as follows:
\begin{eqnarray*}
\varphi_{L,S} (x) \equiv \left \{ \begin{array}{l l}  \frac{x^2}{L^2 S}, & \mbox{if}  \ \ x \le LS, \\
                                                             \frac{x}{L} \left ( \log \left ( \frac{ex}{LS} \right ) \right )^{1/2} \, & \mbox{if} \ \ x \ge LS .
                                              \end{array} \right .
\end{eqnarray*}
Because of the square-root on the log term, this can be regarded as corresponding to a ``sub - Bennett'' type exponential bound.
One of the interesting properties of $\varphi_{L,S}$ established by \cite{MR1258865} is given in the following lemma:
\medskip

\begin{lem}
\label{lem:TalagrandSubBennettFunction}
There is a number $K(L)$ depending on $L$ only such that
\begin{eqnarray*}
\varphi_{L,S} (K(L) x \sqrt{S}) \ge 11 x^2 \ \ \ \mbox{for all} \ \ x \le K(L) \sqrt{S} .
\end{eqnarray*}
\end{lem}
\smallskip

\par\noindent
This is Lemma 3.6 of \cite{MR1258865} page 47.  Talagrand uses this Lemma 
to develop a Kiefer-type inequality:  see also 
\cite{MR1385671}, Corollary A.6.3.
In the basic Kiefer type inequality for Binomial random variables, \cite{MR1385671}, Corollary A.6.3, it follows that
\begin{eqnarray*}
P( \sqrt{n} | \overline{Y}_n - p | \ge z) \le 2\exp(- 11 z^2)  
\end{eqnarray*}
for $\log(1/p) -1 \ge 11$;  i.e. for $p \le e^{-12}$.

A similar fact holds for any exponential bound of the Bennett type
under a certain boundedness hypothesis.  Suppose that 
\begin{eqnarray*}
P( | Z | > z ) \le 2 \exp \left ( - \frac{z^2}{2\tau} \psi \left (\frac{L z}{\tau} \right ) \right) 
\end{eqnarray*} 
and that $P(|Z| \ge v ) = 0$ for all $v \ge C$.  
Then, since $\psi$ is decreasing, for $z \le C$
\begin{eqnarray*}
P( | Z | > z) & \le & 2 \exp \left ( - \frac{z^2}{2\tau} \psi \left ( \frac{LC}{\tau} \right ) \right ) \\
& = & 2 \exp \left ( - \frac{z^2}{2LC} \frac{LC}{\tau} \psi \left ( \frac{LC}{\tau} \right ) \right ) \\
& \le & 2 \exp \left ( - \frac{z^2}{LC}   \log \left ( \frac{LC}{e \tau} \right )  \right ) 
\end{eqnarray*}
where the $\log$ term can be made arbitrarily large by choosing $\tau $  sufficiently small.   
Here the second inequality follows from the fact that 
\begin{eqnarray}
x \psi (x) \ge 2 \log (x/e) = 2 \left ( \log x -1 \right ).
\label{ImprovedPsiLowerBound}
\end{eqnarray}  
\par\noindent
{\bf Proof of (\ref{ImprovedPsiLowerBound}):}  
Since $\psi (x) = 2 x^{-2} h(1+x)$ where $h(x) = x(\log x-1) +1$, we can write, with 
$\underline{\psi} (x) \equiv 2 x^{-2} \log (x/e)$,
\begin{eqnarray*}
\frac{1}{2} \left (x \psi (x) - \underline{\psi} (x) \right ) 
& = & \frac{1}{x} \left \{ (1+x) \log (1+x) -x  \right \} - \left ( \log x -1 \right ) \\
& = & \frac{1+x}{x} \log (1+x) - \log x \\
& = & \log \left ( \frac{1+x}{x} \right ) + \frac{1}{x} \log (1+x)  
\end{eqnarray*}
where both terms are clearly non-negative. 
\hfill $\Box$
\medskip

Now we consider another basic inequality due to \cite{MR1258865}.
Suppose that 
$$
\theta : 2^{\cal X}\cap \{ C\in 2^{\cal X} : \ | C | < \infty \} \mapsto \RR
$$ 
satisfies the following three properties:\\
(a)  $C \subset D $ implies that $\theta (C) \le \theta (D)$ for $C, D \in 2^{\cal X}$.\\
(b)  $\theta (C \cup D) \le \theta (C) + \theta (D)$.  \\
(c)  $\theta (C) \le | C | = \# (C)$.\\
Then if $X_1, \ldots , X_n $ are i.i.d. $P$ non-atomic on $({\cal X}, {\cal A})$ and $Z \equiv \theta ( \{ X_1, \ldots , X_n \})$,
for some universal constant $K_2$ we have, for $z \ge K_2 E(Z)$, 
\begin{eqnarray*}
P( Z \ge z ) \le \exp \left ( - \frac{z}{K_2} \log \left ( \frac{e z}{K_2 E(Z)} \right ) \right ). 
\end{eqnarray*}
As noted by \cite{MR1258865}, this follows from an isoperimetric inequality established in  
\cite{MR1048946}, 
but it is also a consequence of results of \cite{MR1122615,MR1361756}.  
Here we simply note that it can be rephrased as a Bennett type inequality:
for all $z \ge K_2 E(Z)$ 
\begin{eqnarray*}
P( Z \ge z ) \le \exp \left ( - E(Z) h \left ( 1 + \frac{z}{K_2 E(Z)} \right ) \right ) .
\end{eqnarray*}
This follows by simply checking that 
\begin{eqnarray*}
E(Z) h \left ( 1 + \frac{z}{K_2 E(Z)} \right ) \le \frac{z}{K_2 E(Z)} \log \left ( \frac{e z}{K_2 E(Z)} \right ) 
\end{eqnarray*}
for $z \ge K_2 E(Z)$.  
\medskip

Also see
\cite{MR1849347}, Theorem 7.5, page 142 and  Corollary 7.8, page 148; 
\cite{MR1782276}, 
and 
 \cite{MR3185193}, Theorem 6.12, page 182.
\medskip

One further remark seems to be in order:  
\cite{MR1048946} 
Theorem 2 and Proposition 12, shows that Orlicz norms of the Bennett type are ``too large'' to yield nice
generalizations of the classical Hoffmann-J{\o}rgensen inequality in the setting of sums of 
independent bounded sequences in a general Banach space.  This follows by noting that 
Talagrand's condition (2.11) fails for the Bennett-Orlicz norm $\Psi_2 (\cdot, L)$ as defined in (\ref{BennettOrliczNorm}).

\section{Appendix 1:  Lambert's function $W$;  inverses of $h$ and $h_2$}
\label{sec:Appendix1}

Let $h(x) \equiv x (\log x -1)+1$ and $h_2 (x) \equiv h(1+x)$  for $x\ge 0$.  
The function $h$ is convex, decreasing on $[0,1]$, increasing on $[1,\infty)$, with $h(1) = 0$;
see \cite{MR838963}, page 439.  
 The Lambert, or product log function, $W$
(see e.g. \cite{MR1414285} and 
satisfies $W(x) e^{W(x)} = x$ 
for $x \ge -1/e$.  As noted by \cite{MR3185193}, problem 2.18,
the inverse functions $h^{-1} $ (for the function $h: [1,\infty) \rightarrow [0,\infty)$)
and $h_2^{-1} $  (for the function $h_2 : [0,\infty) \rightarrow [0,\infty)$) can be expressed 
in terms of the function $W$.   Here are some facts about $W$:\\
\par\noindent
{\bf Fact 1:}  $W: [-1/e, \infty) \mapsto \RR$ is multi-valued on $[-1/e,0)$ with 
two branches $W_0$ and \\
$\phantom{blablablab}$$W_{-1}$ where $W_0 (x) > 0$, $W_{-1} (x) < 0$, 
and $W_0(-1/e) = -1 = W_{-1} (-1/e)$.\\
{\bf Fact 2:}  $W_0$ is monotone increasing on $[-1/e,\infty)$ with $W(0) = 0$ and $W' (0) =1$.
\smallskip 

\par\noindent
See \cite{MR2655344}, section 4.13, page 111;  
and \cite{MR1414285}.
\medskip

In the following we simply write $W$ for $W_0$.   The following lemma shows that the
inverses of the functions $h$ and $h_2$ can be expressed in terms of $W$.
\medskip
 
\begin{lem}
\label{lem:h2InverseUsingLambertW}
($h$ and $h_2$ inverses in terms of $W$) \\
(i) \ For $y \ge 0$
\begin{eqnarray}
h^{-1} (y) = e \frac{(y-1)/e}{W((y-1)/e)}  \ge 1\; \label{hInverse}
\end{eqnarray}
(ii) \ For $y \ge 0$
\begin{eqnarray}
h_2^{-1} (y) = h^{-1} (y) -1 = e \frac{(y-1)/e}{W((y-1)/e)} -1  \ge 0.  \label{h2Inverse}
\end{eqnarray}
\end{lem}

\par\noindent
{\bf Proof.}   If $h^{-1}$ is as in the display we have, since $h(x) = x(\log x -1 ) +1$,
\begin{eqnarray*}
h( h^{-1} (y)) 
& = & \frac{ e(y-1)/e}{W((y-1)/e)} \left ( \log \left ( \frac{ (y-1)/e}{W((y-1)/e)}\right ) +1 - 1\right ) + 1 \\
& = & \frac{e (y-1)/e}{W((y-1)/e)}  W((y-1)/e)  + 1 
    \ \ \ \ \mbox{since} \ \ e^{W(x)} = \frac{x}{W(x)} \\
& = &  y-1 + 1 = y .
\end{eqnarray*}
Thus (\ref{hInverse}) holds.   Then (\ref{h2Inverse}) follows immediately.  
\hfill $\Box$
\medskip

In view of Lemma~\ref{lem:h2InverseUsingLambertW}, the following lower 
bounds on the function $W$ will be useful in deriving upper 
bounds on $h^{-1}$ and $h_2^{-1}$.
\medskip

\begin{lem}
\label{lem:LowerBoundForW}
(A lower bound for $W$) For $z>0$ 
 \begin{eqnarray}
W(z) \ge \frac{1}{2} \log (e z) = 2^{-1} (1 + \log z) . 
\label{WLowerBoundTwo}
\end{eqnarray}
\end{lem}
\smallskip

\par\noindent
{\bf Proof.}   We first prove (\ref{WLowerBoundTwo}) for $z \ge 1/e$.  
Since $W(z)$ is increasing for $z \ge  0$, the claimed inequality is equivalent to 
\begin{eqnarray*}
z = W(z) e^{W(z)} \ge \frac{1}{2} \log (ez) \exp ( (1/2) \log (ez) ) = (ez)^{1/2} \log ((ez)^{1/2} ) \equiv y \log y,
\end{eqnarray*}
for $ez \ge 1$
where $y \equiv (e z)^{1/2}$.  But then the last display is equivalent to 
$$
e^{-1} y^2 \ge y \log y  \ \ \ \mbox{for} \ \ y \ge 1 
$$
or
$$
g(y) \equiv y^2 - e y \log y \ge 0 \ \ \mbox{for all} \ \ y \ge 1.
$$
Now $g(1) = 0$, $g(e) = 0$, and $g' (y) = 2y - e - e \log y $ has $g'(1) = 2-e<0$, $g' (e) = 0$, 
and $g' (y) > 0$ for $y> e$ with $g'' (y) = 2- e/y$, we find that $g'' (e) = 2 -e/e = 1> 0$.   
Thus the claimed bound holds for $z \ge 1/e$.  For $0 \le z < 1/e$ the bound holds trivially 
since $W(z) \ge 0$ while $2^{-1} \log (ez) < 0$.  
\hfill $\Box$
\medskip

Combining Lemma~\ref{lem:h2InverseUsingLambertW} with the lower bounds for $W$ given in 
Lemma~\ref{lem:LowerBoundForW} yields
the following upper bounds for $h^{-1}$ and $h_2^{-1}$.  The second and third 
parts of the following lemma are motivated by the fact that 
$h_2 (x) = h(1+x) \equiv (x^2/2) \psi (x) $ where $\psi (x) \nearrow 1$ as $x\searrow 0$; 
see \cite{MR838963}, Proposition 4.4.1, page 441.  
\medskip
 
\begin{lem}
\label{lem:UpperBoundsForHInverseandH2Inverse}
(Upper bounds for $h^{-1} $ and $h_2^{-1}$ )\\
(i) \ \ For $y> 1+e$ 
\begin{eqnarray}
h^{-1} (y) \le  2 \frac{y-1}{\log (y-1)} .
\label{h2InverseUpperBoundPoiss}
\end{eqnarray}
(ii) \ For $y > 1+e$,
\begin{eqnarray}
h_2^{-1} (y) \le 2 \frac{y-1}{\log (y-1)} -1 
\label{h2InverseUpperBoundPoiss}
\end{eqnarray}
(iii) For $ 0 \le y \le 9 c^{-2} (c^2/2 -1)^2$ with $c> \sqrt{2}$,  
\begin{eqnarray}
h_2^{-1} (y) \le c \sqrt{y} .
\label{h2InverseUpperBoundGauss}
\end{eqnarray}
In particular, with  $c=2$, the bound holds for $0\le y \le 9/4$, and with $c=2.2$, the bound holds for $0 \le y \le 1+e$.\\
(iv) \ For $0 < y < \infty $,
\begin{eqnarray}
h_2^{-1} (y) \le  \left \{  \begin{array}{l l}  2.2 \sqrt{v} , &  \ \ y \le 1+e, \\ 
                                     2 \frac{y-1}{\log (y-1)} -1 , & \ \ y> 1+e .\end{array} \right .
\label{h2InversePoiss-GaussTrade}
\end{eqnarray}
\end{lem}

\par\noindent
{\bf Proof.}  (i) Follows from (i) of Lemma~\ref{lem:h2InverseUsingLambertW} 
together with Lemma~\ref{lem:LowerBoundForW}. 
Note that 
$g(x) \equiv x/ \log (x) \ge e$ and 
$g$ is increasing for $x \ge e$. \\
(ii) follows from (ii) of Lemma~\ref{lem:h2InverseUsingLambertW}  
and Lemma~\ref{lem:LowerBoundForW}.  \\ 
(iii) To show that (\ref{h2InverseUpperBoundGauss}) holds, note that the inequality is equivalent to 
$y \le h_2 (c\sqrt{y})$, and hence, by taking $x\equiv c \sqrt{y} $, to the inequality 
$$
x^2/c^2 \le h_2 (x) = h(1+x) = \frac{x^2}{2} \psi (x)
$$
where $\psi(x) \equiv (2/x^2) h(1+x) \ge 1/(1+x/3)$ 
by Lemma~\ref{Prohorov-Bennett-Comparison} (iva) (or by (10) of Proposition 11.4.1, 
\cite{MR838963} page 441).  But then we have 
\begin{eqnarray*}
h_2 (x) = h(1+x) \ge \frac{x^2}{2} \frac{1}{1+x/3} \ge \frac{x^2}{c^2} 
\end{eqnarray*}
where the last inequality holds if $0 \le x \le 3(c^2/2-1)$.  Hence the inequality in (iii) holds 
for $0 \le y \equiv x^2/c^2 \le 9 (c^2/2-1)^2/c^2$.  Finally, (iv) holds by combining the bounds in (ii) and (iii).
\hfill $\Box$

\section{Appendix 2:  General versions of Lemmas 1-5}
\label{Appendix2}

Now consider Young functions of the form $\Psi = e^{\psi} -1$ where 
$\psi$ is assumed to be convex and nondecreasing with $\psi (0) = 0$.  
(Note that we have changed notation in this section:  the functions $h$ and $h_j$ for $j\in \{0,1,2,4,6 \}$ in 
Sections~\ref{sec:intro} - \ref{sec:Appendix1} are denoted here by $\psi$.)
Our goal in this section is to give general versions of Lemmas 1 - 5 of \cite{MR3101846}  
and section 3 above.   The advantage of this formulation is that the resulting
lemmas apply to all the special cases treated in Sections 2 and 3 and more.

\begin{lem}
\label{lem:Lem1-genCase}
Suppose that $\tau \equiv \| Z \|_{\Psi} < \infty$.  Then for all $t > 0$ 
\begin{eqnarray*}
P( | Z | > \tau \psi^{-1} (t) ) \le 2 e^{-t} .
\end{eqnarray*}
\end{lem}

For the general version of Lemma 2 we consider a scaled version of $\Psi$ as follows:
\begin{eqnarray}
\Psi (z; L) \equiv \Psi_L (z) \equiv \exp \left ( \frac {2}{L^2} \psi (L z) \right ) -1 .
\label{ScaledPsi}
\end{eqnarray}

\begin{lem}
\label{lem:Lem2-genCase}
Suppose that for some $\tau > 0$ and $L>0$ 
\begin{eqnarray*}
P\left ( | Z| \ge \frac{\tau}{L} \psi^{-1} ( L^2 t/2) \right )  \le 2 e^{-t} \ \ \mbox{for all} \ \ t > 0 .
\end{eqnarray*}
Then  $\| Z \|_{\Psi ( \cdot ; \sqrt{3} L ) } \le \sqrt{3} \tau $.
\end{lem}

\begin{lem}
\label{lem:Lem3-genCase}
Suppose that $\Psi $ is non-decreasing, convex, with $\Psi (0) = 0$. 
Suppose that $Z_1, \ldots , Z_m $ are random variables with $\max_{1 \le j \le m} \| Z_j \|_{\Psi} \equiv \tau < \infty$.
Then
\begin{eqnarray*}
E \left ( \max_{1 \le j \le m} | Z_j | \right ) \le \tau \Psi^{-1} (m) . 
\end{eqnarray*}
\end{lem}

\begin{lem}
\label{lem:Lem4-genCase}
Suppose that $\Psi $ is non-decreasing, convex, with $\Psi (0) = 0$. 
Suppose that $Z_1, \ldots , Z_m $ are random variables with $\max_{1 \le j \le m} \| Z_j \|_{\Psi} \equiv \tau < \infty$.
Then
\begin{eqnarray*}
P \left ( \max_{1 \le j \le m} | Z_j | \ge \tau \left ( \psi^{-1} (\log(1+m)) + \psi^{-1} (t) \right ) \right ) \le 2 e^{-t}  . 
\end{eqnarray*}
\end{lem}

\begin{lem}
\label{lem:Lem5-genCase}
Suppose that $\Psi $ is non-decreasing, convex, with $\Psi (0) = 0$. 
Suppose that $Z_1, \ldots , Z_m $ are random variables with $\max_{1 \le j \le m} \| Z_j \|_{\Psi} \equiv \tau < \infty$.
Then
\begin{eqnarray*}
\bigg \| \left ( \max_{1 \le j \le m} | Z_j | - \tau \Psi^{-1} (m) \right )_+\bigg \|_{\Psi (\cdot ; \sqrt{6})} \le \sqrt{6} \tau . 
\end{eqnarray*}
\end{lem}

\par\noindent
{\bf Proof of Lemma~\ref{lem:Lem1-genCase}.}\ \ 
For all $c> \| Z \|_{\Psi}$
\begin{eqnarray*}
P( | Z | / c > \psi^{-1} (t) ) 
& = & P( \psi ( | Z | / c ) > t ) = P( e^{\psi ( | Z | /c )}-1> e^t -1 ) \\
& = & P( \Psi ( | Z | /c) > e^t -1 ) \le \left ( E \Psi ( | Z | /c ) +1 \right )  e^{-t} .
\end{eqnarray*}
Thus letting $c \searrow \tau$ yields
\begin{eqnarray*}
P( | Z | / \tau > \psi^{-1} (t) ) 
& = & \lim_{c \searrow \tau} P( | Z | /c > \psi^{-1} (t)) \le \lim_{c \searrow \tau} \left ( E \Psi ( | Z |/c) +1 \right ) e^{-t} \le 2 e^{-t} .
\end{eqnarray*}
\hfill $\Box$

\par\noindent
{\bf Proof of Lemma~\ref{lem:Lem2-genCase}.} \ \ 
Let $\alpha, \beta > 0$.   We compute
\begin{eqnarray*}
E \Psi \left ( \frac{|Z|}{\alpha \tau} ; \beta L \right ) 
& = & \int_0^{\infty} P \left ( \Psi \left ( \frac{|Z|}{\alpha \tau}; \beta L \right ) \ge v \right ) dv \\
& = & \int_0^{\infty} P \left ( \frac{2}{\beta^2 L^2} \psi \left ( \frac{\beta L |Z|}{\alpha \tau}\right ) \ge \log (1+v) \right ) dv \\
& = & \int_0^{\infty} P \left ( \frac{\beta L |Z|}{\alpha \tau} \ge \psi^{-1}  \left ( \frac{\beta^2 L ^2 }{2}  \log (1+v) \right ) \right ) dv \\
& = & \int_0^{\infty} P \left ( | Z | \ge \frac{\alpha \tau}{\beta L} \psi^{-1}  \left ( \frac{\beta^2 L ^2 }{2}  \log (1+v) \right ) \right ) dv \\
& = &  \int_0^{\infty} P \left ( | Z | \ge \frac{\alpha \tau}{\beta L} \psi^{-1}  \left ( \frac{\beta^2 L ^2 }{2}  t \right ) \right ) e^t dt  \\
& \le & \int_0^{\infty} 2 e^{- 3 t} e^t dt  = 1
\end{eqnarray*}
by choosing $\alpha = \beta = \sqrt{3} $. 
\hfill $\Box$
\smallskip

\par\noindent
{\bf Proof of Lemma~\ref{lem:Lem3-genCase}.} \ \  Let $c > \tau$.  Then by Jensen's inequality and convexity of $\Psi$
\begin{eqnarray*}
\Psi \left ( \frac{ E \{ \max_{1\le j \le m} | Z_j |}{c} \right ) 
& \le & E \left \{ \Psi\left ( \max_{1 \le j \le m} | Z_j |/c \right ) \right \} = E \left \{ \max_{1 \le j \le m } \Psi ( |Z_j | / c ) \right \} \\
& \le & E \left \{ \sum_{j=1}^m  E \Psi ( | Z_j | / c ) \right \} \le m \cdot \max_{1\le j \le m} E \Psi ( | Z_j | /c ) .
\end{eqnarray*}
Letting $c \searrow \tau$ yields
\begin{eqnarray*}
E \left \{ \max_{1 \le j \le m} | Z_j | \right \} \le \tau \Psi^{-1} \left ( m \cdot \max_{1 \le j \le m}  \Psi ( | Z_j | / \tau ) \right ) = \tau \Psi^{-1} (m ) .
\end{eqnarray*}
\hfill $\Box$
\smallskip

\par\noindent
{\bf Proof of Lemma~\ref{lem:Lem4-genCase}.}  \ \   
For any $u>0$ and $v>0$ concavity of $\psi^{-1}$ implies that 
\begin{eqnarray*}
\psi^{-1} (u) + \psi^{-1} (v) \ge \psi^{-1} (u+v) .
\end{eqnarray*}
Therefore, by using this with $u=\log (1+m)$ and $v= t$,
a union bound, and Lemma~\ref{lem:Lem1-genCase},
\begin{eqnarray*}
\lefteqn{P \left ( \max_{1 \le j \le m} | Z_j | \ge \tau \left ( 
  \psi^{-1} (\log (1+m)) + \psi^{-1} (t) )\right ) \right )} \\
& \le & P \left ( \max_{1 \le j \le m} | Z_j | \ge \tau \psi^{-1} ( \log (1+m) + t ) \right )\\
& \le & \sum_{j=1}^m P\left ( | Z_j | \ge \ \tau \psi^{-1} (   \log (1+m) + t  )\right ) \\
& \le & 2m \exp \left ( - (t + \log (1+m))\right ) = 2 \frac{m}{m+1} e^{-t}  \le 2 e^{-t}.
\end{eqnarray*}
\hfill $\Box$
\smallskip

\par\noindent
{\bf Proof of Lemma~\ref{lem:Lem5-genCase}.}  \ \ 
By Lemma~\ref{lem:Lem4-genCase} 
\begin{eqnarray*}
P \left ( \max_{1 \le j \le m} | Z_j | \ge \tau \left ( 
  \psi^{-1} (\log (1+m)) + \psi^{-1} (t) )\right ) \right ) \le 2 e^{-t}  \ \ \mbox{for all} \ t > 0,
\end{eqnarray*}
so the hypothesis of Lemma~\ref{lem:Lem2-genCase} holds for 
$$
Z \equiv  \left ( \max_{1 \le j \le m} | Z_j | - \tau \psi^{-1} (\log (1+m)) \right )_+
$$
with $L = \sqrt{2}$ and $\tau$ replaced by $\sqrt{2} \tau$.   Thus the conclusion of 
Lemma~\ref{lem:Lem2-genCase} holds for $Z$ with these choices of $L$ and $\tau$:
$\| Z \|_{\Psi (\cdot ; \sqrt{6})}  \le \sqrt{6} \tau $.
\hfill $\Box$
\bigskip

\par\noindent
\textbf{Acknowledgement: }  
I owe thanks to Evan Greene and Johannes Lederer for several helpful conversations 
and suggestions.  
Thanks are also due to Richard Nickl for a query concerning Prokhorov's inequality.

\def\cprime{$'$}


\begin{thebibliography}{}

\bibitem[Arcones and Gin{\'e}(1995)Arcones and Gin{\'e}]{MR1348376}
Arcones, M.~A. and Gin{\'e}, E. (1995).
\newblock On the law of the iterated logarithm for canonical {$U$}-statistics
  and processes.
\newblock {\em Stochastic Process. Appl.}, {\bf 58}(2), 217--245.

\bibitem[Bennett(1962)Bennett]{Bennett:62}
Bennett, G. (1962).
\newblock Probability inequalities for the sum of independent random variables.
\newblock {\em Journal of the American Statistical Association\/}, {\bf 57},
  33--45.

\bibitem[Birg{\'e} and Massart(1998)Birg{\'e} and Massart]{MR1653272}
Birg{\'e}, L. and Massart, P. (1998).
\newblock Minimum contrast estimators on sieves: exponential bounds and rates
  of convergence.
\newblock {\em Bernoulli\/}, {\bf 4}(3), 329--375.

\bibitem[Boucheron {\em et~al.}(2013)Boucheron, Lugosi, and Massart]{MR3185193}
Boucheron, S., Lugosi, G., and Massart, P. (2013).
\newblock {\em Concentration Inequalities\/}.
\newblock Oxford University Press, Oxford.

\bibitem[Corless {\em et~al.}(1996)Corless, Gonnet, Hare, Jeffrey, and
  Knuth]{MR1414285}
Corless, R.~M., Gonnet, G.~H., Hare, D. E.~G., Jeffrey, D.~J., and Knuth, D.~E.
  (1996).
\newblock On the {L}ambert {$W$} function.
\newblock {\em Adv. Comput. Math.}, {\bf 5}(4), 329--359.

\bibitem[de~la Pe{\~n}a and Gin{\'e}(1999)de~la Pe{\~n}a and
  Gin{\'e}]{MR1666908}
de~la Pe{\~n}a, V.~H. and Gin{\'e}, E. (1999).
\newblock {\em Decoupling; From dependence to independence\/}.
\newblock Probability and its Applications (New York). Springer-Verlag, New
  York.

\bibitem[Dudley(1999)Dudley]{MR1720712}
Dudley, R.~M. (1999).
\newblock {\em Uniform Central Limit Theorems\/}, volume~63 of {\em Cambridge
  Studies in Advanced Mathematics\/}.
\newblock Cambridge University Press, Cambridge.

\bibitem[Ghosh and Goldstein(2011a)Ghosh and Goldstein]{MR2763529}
Ghosh, S. and Goldstein, L. (2011a).
\newblock Applications of size biased couplings for concentration of measures.
\newblock {\em Electron. Commun. Probab.}, {\bf 16}, 70--83.

\bibitem[Ghosh and Goldstein(2011b)Ghosh and Goldstein]{MR2773032}
Ghosh, S. and Goldstein, L. (2011b).
\newblock Concentration of measures via size-biased couplings.
\newblock {\em Probab. Theory Related Fields\/}, {\bf 149}(1-2), 271--278.

\bibitem[Goldstein and I{\c{s}}lak(2014)Goldstein and I{\c{s}}lak]{MR3162712}
Goldstein, L. and I{\c{s}}lak, {\"U}. (2014).
\newblock Concentration inequalities via zero bias couplings.
\newblock {\em Statist. Probab. Lett.}, {\bf 86}, 17--23.

\bibitem[Hewitt and Stromberg(1975)Hewitt and Stromberg]{MR0367121}
Hewitt, E. and Stromberg, K. (1975).
\newblock {\em Real and Abstract Analysis\/}.
\newblock Springer-Verlag, New York-Heidelberg.
\newblock A modern treatment of the theory of functions of a real variable,
  Third printing, Graduate Texts in Mathematics, No. 25.

\bibitem[Johnson {\em et~al.}(1985)Johnson, Schechtman, and Zinn]{MR770640}
Johnson, W.~B., Schechtman, G., and Zinn, J. (1985).
\newblock Best constants in moment inequalities for linear combinations of
  independent and exchangeable random variables.
\newblock {\em Ann. Probab.}, {\bf 13}(1), 234--253.

\bibitem[Krasnosel{\cprime}ski{\u\i} and
  Ruticki{\u\i}(1961)Krasnosel{\cprime}ski{\u\i} and Ruticki{\u\i}]{MR0126722}
Krasnosel{\cprime}ski{\u\i}, M.~A. and Ruticki{\u\i}, J.~B. (1961).
\newblock {\em Convex Functions and {O}rlicz Spaces\/}.
\newblock Translated from the first Russian edition by Leo F. Boron. P.
  Noordhoff Ltd., Groningen.

\bibitem[Kruglov(2006)Kruglov]{MR2331988t}
Kruglov, V.~M. (2006).
\newblock Strengthening of {P}rokhorov's arcsine inequality.
\newblock {\em Theor. Probab. Appl.}, {\bf 50}, 677 -- 684.
\newblock Transl. from Strengthening the {P}rokhorov arcsine inequality, {\sl
  Teor. Veroyatn. Primen.}, {\bf 50}, (2005).

\bibitem[Ledoux(2001)Ledoux]{MR1849347}
Ledoux, M. (2001).
\newblock {\em The Concentration of Measure Phenomenon\/}, volume~89 of {\em
  Mathematical Surveys and Monographs\/}.
\newblock American Mathematical Society, Providence, RI.

\bibitem[Massart(2000)Massart]{MR1782276}
Massart, P. (2000).
\newblock About the constants in {T}alagrand's concentration inequalities for
  empirical processes.
\newblock {\em Ann. Probab.}, {\bf 28}(2), 863--884.

\bibitem[Pisier(1983)Pisier]{MR717231}
Pisier, G. (1983).
\newblock Some applications of the metric entropy condition to harmonic
  analysis.
\newblock In {\em Banach spaces, harmonic analysis, and probability theory
  ({S}torrs, {C}onn., 1980/1981)\/}, volume 995 of {\em Lecture Notes in
  Math.}, pages 123--154. Springer, Berlin.

\bibitem[Pollard(1990)Pollard]{MR1089429}
Pollard, D. (1990).
\newblock {\em Empirical Processes: Theory and Applications\/}.
\newblock NSF-CBMS Regional Conference Series in Probability and Statistics, 2.
  Institute of Mathematical Statistics, Hayward, CA; American Statistical
  Association, Alexandria, VA.

\bibitem[Prokhorov(1959)Prokhorov]{MR0121857}
Prokhorov, Y.~V. (1959).
\newblock An extremal problem in probability theory.
\newblock {\em Theor. Probability Appl.}, {\bf 4}, 201--203.

\bibitem[Rao and Ren(1991)Rao and Ren]{MR1113700}
Rao, M.~M. and Ren, Z.~D. (1991).
\newblock {\em Theory of {O}rlicz spaces\/}, volume 146 of {\em Monographs and
  Textbooks in Pure and Applied Mathematics\/}.
\newblock Marcel Dekker, Inc., New York.

\bibitem[Roy and Olver(2010)Roy and Olver]{MR2655344}
Roy, R. and Olver, F. W.~J. (2010).
\newblock Elementary functions.
\newblock In {\em N{IST} handbook of mathematical functions\/}, pages 103--134.
  U.S. Dept. Commerce, Washington, DC.

\bibitem[Shorack and Wellner(1986)Shorack and Wellner]{MR838963}
Shorack, G.~R. and Wellner, J.~A. (1986).
\newblock {\em Empirical {P}rocesses with {A}pplications to {S}tatistics\/}.
\newblock Wiley Series in Probability and Mathematical Statistics: Probability
  and Mathematical Statistics. John Wiley \& Sons Inc., New York.

\bibitem[Stout(1974)Stout]{MR0455094}
Stout, W.~F. (1974).
\newblock {\em Almost {S}ure {C}onvergence\/}.
\newblock Academic Press [A subsidiary of Harcourt Brace Jovanovich,
  Publishers], New York-London.
\newblock Probability and Mathematical Statistics, Vol. 24.

\bibitem[Talagrand(1989)Talagrand]{MR1048946}
Talagrand, M. (1989).
\newblock Isoperimetry and integrability of the sum of independent
  {B}anach-space valued random variables.
\newblock {\em Ann. Probab.}, {\bf 17}(4), 1546--1570.

\bibitem[Talagrand(1991)Talagrand]{MR1122615}
Talagrand, M. (1991).
\newblock A new isoperimetric inequality and the concentration of measure
  phenomenon.
\newblock In {\em Geometric aspects of functional analysis (1989--90)\/},
  volume 1469 of {\em Lecture Notes in Math.}, pages 94--124. Springer, Berlin.

\bibitem[Talagrand(1994)Talagrand]{MR1258865}
Talagrand, M. (1994).
\newblock Sharper bounds for {G}aussian and empirical processes.
\newblock {\em Ann. Probab.}, {\bf 22}(1), 28--76.

\bibitem[Talagrand(1995)Talagrand]{MR1361756}
Talagrand, M. (1995).
\newblock Concentration of measure and isoperimetric inequalities in product
  spaces.
\newblock {\em Inst. Hautes \'Etudes Sci. Publ. Math.}, {\bf 81}, 73--205.

\bibitem[van~de Geer and Lederer(2013)van~de Geer and Lederer]{MR3101846}
van~de Geer, S. and Lederer, J. (2013).
\newblock The {B}ernstein-{O}rlicz norm and deviation inequalities.
\newblock {\em Probab. Theory Related Fields\/}, {\bf 157}(1-2), 225--250.

\bibitem[van~der Vaart and Wellner(1996)van~der Vaart and Wellner]{MR1385671}
van~der Vaart, A.~W. and Wellner, J.~A. (1996).
\newblock {\em Weak {C}onvergence and {E}mpirical {P}rocesses\/}.
\newblock Springer Series in Statistics. Springer-Verlag, New York.

\end{thebibliography}

\end{document}